\documentclass[12pt]{article} 
\usepackage{amssymb} 
\usepackage{amsfonts} 
\usepackage{amsmath} 
\input{epsf.sty} 
\usepackage[dvips]{graphicx} 
\usepackage{latexsym,amsmath,amsfonts,amscd} 
\usepackage{epsfig} 
\usepackage{changebar} 
\usepackage{fancybox} 
\title 
{Gradient Estimates for 
 the Perfect Conductivity Problem} 
\author{Ellen Shiting Bao\\ 
YanYan Li 
\thanks{Partially supported by NSF grant DMS-0401118.} 
\\Biao Yin 
    \\ Department of Mathematics 
\\ 
  Rutgers University\\ 
     110 Frelinghuysen Rd. \\ 
        Piscataway, NJ 08854} 
\date{} 
 
\begin{document} 
\newtheorem{Def}{Definition}[section] 
\newtheorem{thm}{Theorem}[section] 
\newtheorem{lem}{Lemma}[section] 
\newtheorem{rem}{Remark}[section] 
\newtheorem{prop}{Proposition}[section] 
\newtheorem{cor}{Corollary}[section] 
\def\av{{\int \hspace{-2.25ex}-} } 
\maketitle 
 
\bigskip 
 
\setcounter{section}{-1} 
%%%%%%%%%%%%%%%%%%%%%%%%%%%%%%%%%%%%%%%%%%%%%%%%%%%%%%%%%%%%%%%%%%%%%%%% 
% 
%     1. Introduction and Statements of Results 
% 
%%%%%%%%%%%%%%%%%%%%%%%%%%%%%%%%%%%%%%%%%%%%%%%%%%%%%%%%%%%%%%%%%%%%%%%% 
 
\section{Introduction} 
Let $\Omega$ be a bounded open set in $\mathbb{R}^n$ with 
$C^{2,\alpha}$ boundary, $n\geq 2$, $0<\alpha <1$, $D_{1}$ and 
$D_{2}$ be two bounded strictly convex open subsets in $\Omega$ with 
$C^{2,\alpha}$ boundaries which are $\varepsilon$ apart and far away 
from $\partial\Omega$, i.e. 
\begin{equation} \label{domain conditions} 
\begin{split} 
&~~~~\overline D_1, \overline D_2\subset \Omega, ~~~\text{the 
principal 
 curvature of} ~\partial D_1, \partial D_2\geq \kappa_0\\ 
&\varepsilon:=\text{dist}(D_1, D_2)>0, ~~~\text{dist} (D_1\cup D_2, 
\partial\Omega)>r_0,~~~ \text{diam}(\Omega)<\frac 1{r_0}, 
 \end{split} 
\end{equation} 
where $ \kappa_0, 
r_0>0$ are universal constants independent of $\varepsilon$.\\ 
We denote 
$$ 
\widetilde\Omega :=\Omega\backslash\overline{D_{1}\cup D_{2}}. 
$$ 
Given $\varphi\in C^{2}(\partial\Omega)$, consider the following 
scalar equation with Dirichlet boundary condition: 
\begin{equation} \label{eq:k finite} 
\left\{ \begin{aligned} 
 &div(a_k(x)\nabla u_k)=0 \hspace{1cm}in\hspace{0.3cm}\Omega, \\ 
           &u_k=\varphi 
           \hspace{3.1cm}on\hspace{0.2cm}\partial\Omega, 
                 \end{aligned} 
\right. 
\end{equation} 
 where 
\begin{equation} 
            a_k(x)=\left\{ \begin{aligned} 
           &k\in (0, \infty)\hspace{1cm}in\hspace{0.3cm}D_{1} \cup D_{2}, \\ 
           &1\hspace{1cm}on\hspace{0.3cm}\Omega\backslash\overline{D_{1}\cup 
           D_{2}}. 
          \end{aligned} 
\right. \label{0.3} 
\end{equation} 
\\ 
It is well known that there exists a unique solution $u_k \in 
H^1(\Omega)$ of the above equation, which is also the minimizer of 
$I_{k}$ on $H^{1}_{\varphi}(\Omega)$, where 
$$ 
H^{1}_\varphi(\Omega):=\{u\in H^{1}(\Omega)~|~u=\varphi ~on~ 
\partial\Omega\}, 
\qquad I_{k}[v]:=\frac{1}{2}\int_\Omega a_k|\nabla v|^2.$$ 
 
As explained in the introduction of $\cite{LV}$, the above equation in 
dimension $n=2$ can be used as a simple model in the study of 
composite media with closely spaced interfacial boundaries. For this 
purpose, the domain $\Omega$ would model the cross-section of a 
fiber-reinforced composite,  $D_1$ and $D_2$ would represent the 
cross-sections of the fibers,  $\widetilde \Omega$ would represent 
the matrix surrounding the fibers, and the shear modulus of the 
fibers would be $k$ and that of the matrix would be $1$. Equation 
(\ref{eq:k finite}) is then obtained by using a standard model of 
anti-plane shear, and the solution $u_k$ represents the out of plane 
elastic displacement.  The most important quantities from an 
engineering point of view are the stresses, in this case represented 
by $\nabla u_k$.

It is well known that the solution $u_k$ satisfies $\|u_k\|_{ 
C^{2,\alpha}(D_i)}<\infty$.  In fact, if $\partial D_1$ and 
$\partial D_2$ are $C^{m, \alpha}$,  we have $\|u_k\|_{ 
C^{m,\alpha}(D_i)}<\infty$.  Such results do not require $D_i$ to be 
convex and hold for general elliptic systems with piecewise smooth 
coefficients; see e.g. theorem 9.1 in \cite{LV} and proposition 1.6 
in \cite{LN}.   For a fixed $0<k<\infty$, the $C^{m, 
\alpha}(D_i)$-norm of the solution might tend to infinity as 
$\varepsilon\to 0$.   Babuska, Anderson, Smith and Levin \cite{BASL} 
were interested in linear elliptic systems of elasticity arising 
from the study of composite material. They observed numerically 
that, for solution $u$ to certain homogeneous isotropic linear 
systems of elasticity, $\|\nabla u\|_{ L^\infty }$ is bounded 
independently of the distance $\varepsilon$ between $D_1$ and $D_1$. 
Bonnetier and Vogelius \cite{BV} proved this in dimension $n=2$ for 
the solution $u_k$ of (\ref{eq:k finite}) when $D_1$ and $D_2$ are 
two unit balls touching at a point. This result was extended   by Li 
and Vogelius in \cite{LV} 
 to general second order elliptic 
equations with piecewise smooth coefficients, where stronger $C^{1, 
\beta}$ estimates were established.   The $C^{1, \beta}$ estimates 
were further extended by Li and Nirenberg  in \cite{LN} to general 
second order elliptic systems including systems of elasticity.   For 
higher derivative estimates, e.g. an $\varepsilon$-independent 
$L^\infty$-estimate of second derivatives of $u_k$ in $D_1$, we draw 
attention of readers to the 
 open problem on page 894 of \cite{LN}. In \cite{LV} and \cite{LN}, 
the ellipticity constants are assumed to be away from $0$ and 
$\infty$.  If we allow ellipticity constants to deteriorate, the 
situation is different.  It has been shown in various papers, see 
e.g. \cite{BC} and \cite{M}, that when $k=\infty$ the 
$L^\infty$-norm of $\nabla u_k$ for  the solution $u_k$ of equation 
(\ref{eq:k finite}) generally becomes unbounded as $\varepsilon$ 
tends to zero.  The rate at which the $L^\infty$ norm of the 
gradient of a special solution has been shown in \cite{BC} to be 
$\varepsilon^{-1/2}$. 
 
In this paper, we consider the perfect conductivity problem, where 
$k=+\infty$.   It was proved by Ammari, Kang and Lim in \cite{AKL} 
and Ammari, Kang, H. Lee, J. Lee and Lim in \cite{AKLLL} that, when 
$D_1$ and $D_2$ are balls of comparable radii  embedded in 
$\Omega=\mathbb{R}^2$, the blow-up rate of the gradient of the 
solution to the perfect conductivity problem is $\varepsilon^{-1/2}$ 
as $\varepsilon$ goes to zero; with the lower bound given in 
\cite{AKL} and the upper bound given in \cite{AKLLL}. Yun in 
\cite{Yu} generalized the above mentioned result in \cite{AKL} by 
establishing the same lower bound, $\varepsilon^{-1/2}$, 
 for two strictly convex subdomains in $\mathbb{R}^2$. In 
this paper, we give both lower and upper bounds to blow-up rate of 
the gradient for the solution to the perfect conductivity problem in 
a bounded matrix, where two strictly convex subdomains are embedded. 
Our methods apply to dimension $n\ge 3$ as well.  One might 
reasonably suspect that the blow-up rate in dimension $n\ge 3$ 
should be smaller than that in dimension $n=2$.  However we prove 
 the opposite:  As $\varepsilon$ goes to zero, the blow-up 
rate is $\varepsilon^{-1/2}$, $(\varepsilon|\ln{\varepsilon}|)^{-1}$ 
and $\varepsilon^{-1}$ for $n=2,~3$ and $n\geq 4$, respectively. We 
also give a criteria, in terms of a linear functional of the 
boundary data $\varphi$, for the situation where the rate of blow-up 
is realized. Note that \cite{AKL} and \cite{AKLLL} contain also results for $k<\infty$.\\ 
%%%%%%%%%%%%%%%%%%%%%%%%%%%%%%%%%%%%%%%%%%%%%%%%%%%%%%%%%%%%%%%%%%%%%%%% 
% 
%      Describe the problem 
% 
%%%%%%%%%%%%%%%%%%%%%%%%%%%%%%%%%%%%%%%%%%%%%%%%%%%%%%%%%%%%%%%%%%%%%%%% 
 
The perfect conductivity problem is described as follows: 
\begin{equation} \label{eq:k +infty} 
\left\{ \begin{aligned} 
           &\Delta u=0 \hspace{2.28cm}in\hspace{0.3cm}\widetilde\Omega, \\ 
           &u|_{+}=u|_{-}\hspace{1.85cm}on\hspace{0.2cm}\partial D_1\cup\partial D_2, \\ 
           &\nabla u\equiv 0\hspace{2.22cm}in\hspace{0.2cm} D_1\cup D_2, \\ 
           &\int_{\partial D_i}\frac{\partial u}{\partial\nu}\Big|_{+}=0 
           \hspace{0.94cm}(i=1,2), \\ 
           &u=\varphi\hspace{2.49cm}on\hspace{0.2cm}\partial\Omega. 
          \end{aligned} 
\right. 
%\nonumber 
\end{equation} 
where 
$$ 
\frac{\partial u}{\partial\nu}\Big|_{+}:=\lim_{t\rightarrow 
0^{+}}\frac{u(x+t\nu)-u(x)}{t}. 
$$ 
Here and throughout this paper $\nu$ is the outward unit normal to 
the domain and the subscript $\pm$ indicates the limit from outside 
and inside the domain, 
respectively.\\ 
 
The existence and uniqueness of solutions to equation (\ref{eq:k 
+infty}) are well known, see the Appendix. Moreover, the solution 
$u\in H^1(\Omega)$ is the weak limit of the solutions $u_k$ to 
equations (\ref{eq:k finite}) as $k\rightarrow +\infty$. It can be 
also described as the unique function which has the ``~least energy" 
in appropriate functional space, defined as 
$I_{\infty}[u]=\min_{v\in\mathcal{A}}I_{\infty}[v],$ where 
$$ 
I_{\infty}[v]:=\frac{1}{2}\int_{\widetilde{\Omega}}|\nabla v|^2, 
\hspace{2cm}v\in\mathcal{A}, 
$$ 
$$ 
\mathcal{A}:=\big\{v\in H^1_\varphi(\Omega)\big| \nabla v\equiv 0 
~in~ D_{1}\cup D_{2}\big\}. 
$$ 
The readers can refer to the Appendix for the proofs of the above 
statements.\\ 
%%%%%%%%%%%%%%%%%%%%%%%%%%%%%%%%%%%%%%%%%%%%%%%%%%%%%%%%%%%%%%%%%%%%%%%% 
% 
%      Statements of Results 
% 
%%%%%%%%%%%%%%%%%%%%%%%%%%%%%%%%%%%%%%%%%%%%%%%%%%%%%%%%%%%%%%%%%%%%%%%% 
 
We now state more precisely what it means by saying that the 
boundary of a domain, say $\Omega$,  is $C^{2, \alpha}$ for 
$0<\alpha<1$: In a neighborhood of every point of $\partial \Omega$, 
$\partial \Omega$ is the graph of some  $C^{2, \alpha}$ functions of 
$n-1$ variables.  We define the $C^{2,\alpha}$ norm of $\partial 
\Omega$, denoted as $\|\partial \Omega\|_{C^{2,\alpha}}$, 
 as the smallest positive number $\frac 1a$ such that in the $2a-$neighborhood 
of every point of $\partial \Omega$, identified as $0$ after a 
possible translation and rotation of the coordinates so that $x_n=0$ 
is the tangent  to $\partial \Omega$ at $0$, $\partial \Omega$ is 
given by the graph of a $C^{2,\alpha}$ function, denoted as $f$, 
which is defined as $|x'|<a$, the $a-$neighborhood of $0$ in the 
tangent plane. Moreover, $\|f\|_{  C^{2,\alpha}(|x'|<a)} \le \frac 
1a$. 
 
\begin{thm} \label{thm:upbdd} 
Let $\Omega, D_1, D_2\subset \mathbb{R}^n$, $\varepsilon$ be defined 
as in (\ref{domain conditions}), $\varphi\in C^{2}(\partial\Omega)$. 
Let $u\in H^1(\Omega)\cap C^1(\overline{\widetilde{\Omega}})$ be the 
solution to equation (\ref{eq:k +infty}). For $\varepsilon$ 
sufficiently small, there is a positive constant $C$ which depends 
only 
 on $n$, $\kappa_0$, 
$r_0$,  $\|\partial \Omega\|_{C^{2,\alpha}}$, $\|\partial 
D_1\|_{C^{2,\alpha}}$ and $\|\partial D_2\|_{C^{2,\alpha}}$, but 
independent of $\varepsilon$ such that 
\begin{equation} \label{eq:upper} 
\begin{aligned} 
           &\|\nabla u\|_{L^{\infty}(\widetilde{\Omega})}\leq 
\frac{C}{\sqrt{\varepsilon}}\|\varphi\|_{C^2(\partial \Omega)}~~~~~~~~~for \hspace{0.2cm} n=2,\\ 
           &\|\nabla u\|_{L^{\infty}(\widetilde{\Omega})}\leq 
           \frac{C}{\varepsilon|\ln{\varepsilon}|}\|\varphi\|_{C^2(\partial \Omega)}~~~~~for \hspace{0.2cm} n=3,\\ 
           &\|\nabla u\|_{L^{\infty}(\widetilde{\Omega})}\leq 
           \frac{C}{\varepsilon}\|\varphi\|_{C^2(\partial \Omega)}~~~~~~~~~~for \hspace{0.2cm} n \geq 
           4.\\ 
\end{aligned} 
\end{equation} 
\end{thm} 
\begin{rem}\label{rem1} 
We draw attention of readers to the independent work of Yun 
\cite{Yu2} where he has also established the upper bound, 
$\varepsilon^{-1/2}$, in  $\mathbb{R}^2$. The methods are very 
different.  Results in this paper and those in \cite{Yu} and 
\cite{Yu2} do not really need $D_1$ and $D_2$ to be strictly convex, 
the strict convexity is only needed for the portions in a fixed 
neighborhood (the size of the neighborhood is indepedent of 
$\varepsilon$) of a pair of points on $\partial D_1$ and $\partial 
D_2$ which realize minimal distance $\varepsilon$. 
In fact, our proofs of Theorem \ref{thm:upbdd}$-$\ref{thm:lowbdd} 
also apply, with minor modification, to more general situations where
 two 
inclusions, $D_1$ and $D_2$,
 are not necessarily convex near points 
on the boundaries where minimal distance $\varepsilon$
is realized; see discussions after the proofs of
Theorem \ref{thm:upbdd}-\ref{thm:lowbdd} in
Section 1.3.
\end{rem} 
 
To prove Theorem \ref{thm:upbdd}, we first decompose the solution 
$u$ of equation (\ref{eq:k +infty}) as follows: 
\begin{equation} 
\label{eq:decomp} u=C_1v_1+C_2v_2+v_3 
\end{equation} 
where $C_i:=C_i(\varepsilon)~(i=1,2)$ be the boundary value of $u$ 
on $\partial D_i ~(i=1,2)$ respectively, and $v_i\in 
C^2(\overline{\widetilde\Omega})$ $(i=1,2,3)$ satisfies 
\begin{equation} 
\label{eq:v1} \left\{ \begin{aligned} 
           &\Delta v_1=0 \hspace{2.1cm}in~\widetilde\Omega, \\ 
           &v_1=1~~on~\partial D_1,~~~v_1=0~~on~\partial D_2 \cup \partial\Omega, 
          \end{aligned} 
\right. 
%\nonumber 
\end{equation} 
\begin{equation} 
\label{eq:v2} \left\{ \begin{aligned} 
           &\Delta v_2=0 \hspace{2.1cm}in~\widetilde\Omega, \\ 
           &v_2=1~~on~\partial D_2,~~~v_2=0~~on~\partial D_1 \cup \partial\Omega, 
          \end{aligned} 
\right. 
%\nonumber 
\end{equation} 
\begin{equation} 
\label{eq:v3} \left\{ \begin{aligned} 
           &\Delta v_3=0 \hspace{2.1cm}in~\widetilde\Omega, \\ 
           &v_3=0~~on~\partial D_1\cup\partial 
           D_2,~~~v_3=\varphi~~on~\partial\Omega. 
          \end{aligned} 
\right. 
%\nonumber 
\end{equation} 
Define 
\begin{equation} \label{eq:Q_e} 
Q_{\varepsilon}[\varphi]:=\int_{\partial D_1}\frac{\partial 
v_3}{\partial \nu}\int_{\partial \Omega}\frac{\partial v_2}{\partial 
\nu}-\int_{\partial D_2}\frac{\partial v_3}{\partial 
\nu}\int_{\partial \Omega}\frac{\partial v_1}{\partial \nu}, 
\end{equation} 
then $Q_{\varepsilon}: C^2(\partial\Omega)\rightarrow \mathbb{R}$ is 
a linear functional. 
\begin{thm} \label{thm:lowbdd} 
With the same conditions in Theorem \ref{thm:upbdd}, let $u\in 
H^1(\Omega)\cap C^1(\overline{\widetilde{\Omega}})$ be the solution 
to equation (\ref{eq:k +infty}). For $\varepsilon$ sufficiently 
small, there exists a positive constant $C$ which depends on $n$, 
$\kappa_0$, $r_0$, $\|\partial \Omega\|_{C^{2,\alpha}}$, $\|\partial 
D_1\|_{C^{2,\alpha}}$, $\|\partial D_2\|_{C^{2,\alpha}}$ and 
$\|\varphi\|_{C^2(\partial \Omega)}$, but is independent of 
$\varepsilon$ such that 
\begin{equation} \label{eq:lower} 
\begin{aligned} 
           &\|\nabla u\|_{L^{\infty}(\widetilde{\Omega})}\geq 
\frac{|Q_{\varepsilon}[\varphi]|}{C}\cdot\frac{1}{\sqrt{\varepsilon}}~~~~~~~~for \hspace{0.2cm} n=2,\\ 
           &\|\nabla u\|_{L^{\infty}(\widetilde{\Omega})}\geq 
           \frac{|Q_{\varepsilon}[\varphi]|}{C}\cdot\frac{1}{\varepsilon|\ln{\varepsilon}|}~~~~for \hspace{0.2cm} n=3,\\ 
           &\|\nabla u\|_{L^{\infty}(\widetilde{\Omega})}\geq 
           \frac{|Q_{\varepsilon}[\varphi]|}{C}\cdot\frac{1}{\varepsilon}~~~~~~~~~~for \hspace{0.2cm} n \geq 4.\\ 
\end{aligned} 
\end{equation} 
\end{thm} 
 
\begin{rem} \label{rem:varphi=1} 
If $\varphi\equiv0$, then the solution to equation (\ref{eq:k 
+infty}) is $u\equiv0$. Theorem \ref{thm:upbdd} and Theorem 
\ref{thm:lowbdd} are obvious in this case. So we only need to prove 
them for $\|\varphi\|_{C^2(\partial \Omega)}=1$, by considering 
$u/\|\varphi\|_{C^2(\partial \Omega)}$. 
\end{rem} 
\begin{rem} 
It is interesting to know when 
$|Q_\varepsilon[\varphi]|\geq\frac{1}{C}$ for some positive constant 
$C$ independent of $\varepsilon$. Roughly speaking 
$Q_\varepsilon[\varphi]\rightarrow Q^*[\varphi]$ as $\varepsilon\to 
0$, and this amounts to $Q^*[\varphi]\neq 0$. For details, see 
Section 2. 
\end{rem}

Theorem \ref{thm:upbdd}$-$\ref{thm:lowbdd} can be extended to 
equations with more general coefficients as follows: Let $n$, 
$\Omega$, $D_1$, $D_2$, $\varepsilon$ and $\varphi$ be same as in 
Theorem \ref{thm:upbdd}, and let 
$$A_2(x):=\big(a_2^{ij}(x)\big)\in C^2(\overline{\widetilde\Omega})$$ 
be $n\times n$ symmetric matrix functions in $\widetilde\Omega$ 
satisfying for some constants $0<\lambda\leq\Lambda<\infty$, 
$$~\lambda|\xi|^2\leq a_2^{ij}(x)\xi_i\xi_j\leq \Lambda|\xi|^2, 
~~~~~\forall x\in\widetilde\Omega, ~\forall \xi\in\mathbb{R}^n,$$ 
and $a_2^{ij}(x) \in C^2(\overline{\Omega\backslash\omega})$.\\ 
 
We consider 
\begin{equation} \label{eq:k +infty general} 
\left\{ \begin{aligned} 
           &\partial_{x_j}\Big(a_2^{ij}(x)~\partial_{x_i}{u}\Big)=0~~~~~~~~in~\widetilde\Omega,\\ 
           &u|_{+}=u|_{-}~~~~~~~~~~~~~~~~~~~~~~on~\partial D_1\cup\partial D_2, \\ 
           &\nabla u=0~~~~~~~~~~~~~~~~~~~~~~~~~in~D_1\cup D_2, \\ 
           &\int_{\partial D_i}a_2^{ij}(x)\partial_{x_i}{u}\nu_j\big|_{+}=0 ~~~~(i=1,2),\\ 
           &u=\varphi~~~~~~~~~~~~~~~~~~~~~~~~~~~on~\partial\Omega. 
          \end{aligned} 
\right. 
\end{equation} 
where repeated indices denote as usual summations.\\ 
 
Here is an extension of Theorem \ref{thm:upbdd}: 
\begin{thm}\label{thm:upbdd general} 
With the above assumptions, let $u\in H^1(\Omega)\cap 
C^1(\overline{\widetilde\Omega})$ be the solution to equation 
(\ref{eq:k +infty general}). For $\varepsilon$ sufficient small, 
there is a positive constant $C$ which depends only on $n$, 
$\kappa_0$, $r_0$, $\|\partial \Omega\|_{C^{2,\alpha}}$, $\|\partial 
D_1\|_{C^{2,\alpha}}$, $\|\partial D_2\|_{C^{2,\alpha}}$, $\lambda$, 
$\Lambda$ and $\|A_2\|_{C^2(\overline{\widetilde\Omega})}$, but 
independent of $\varepsilon$ such that estimate (\ref{eq:upper}) 
holds. 
\end{thm} 
 
Similar to the decomposition formula (\ref{eq:decomp}), we decompose 
the solution $u$ of equation (\ref{eq:k +infty general}) as follows: 
\begin{equation} 
\label{eq:decomp general} u=C_1V_1+C_2V_2+V_3 
\end{equation} 
where $C_i:=C_i(\varepsilon)~(i=1,2)$ be the boundary value of $u$ 
on $\partial D_i ~(i=1,2)$ respectively, and $V_i\in 
C^2(\overline{\widetilde\Omega})$ $(i=1,2,3)$ satisfies 
\begin{equation} 
\label{eq:V1 general} \left\{ \begin{aligned} 
           &\partial_{x_j}\Big(a_2^{ij}(x)~\partial_{x_i}{V_1}\Big)=0~~~~~~~~in~\widetilde\Omega,\\ 
           &V_1=1~~on~\partial D_1,~~~V_1=0~~on~\partial D_2 \cup \partial\Omega, 
          \end{aligned} 
\right. 
%\nonumber 
\end{equation} 
\begin{equation} 
\label{eq:V2 general} \left\{ \begin{aligned} 
           &\partial_{x_j}\Big(a_2^{ij}(x)~\partial_{x_i}{V_2}\Big)=0~~~~~~~~in~\widetilde\Omega,\\ 
           &V_2=1~~on~\partial D_2,~~~V_2=0~~on~\partial D_1 \cup \partial\Omega, 
          \end{aligned} 
\right. 
%\nonumber 
\end{equation} 
\begin{equation} 
\label{eq:V3 general} \left\{ \begin{aligned} 
           &\partial_{x_j}\Big(a_2^{ij}(x)~\partial_{x_i}{V_3}\Big)=0~~~~~~~~in~\widetilde\Omega,\\ 
           &V_3=0~~on~\partial D_1\cup\partial 
           D_2,~~~V_3=\varphi~~on~\partial\Omega. 
          \end{aligned} 
\right. 
%\nonumber 
\end{equation} 
Define 
\begin{equation} \label{eq:Q_e general} 
\begin{split} 
Q_{\varepsilon}[\varphi]:=&\int_{\partial 
D_1}a_2^{ij}(x)~\partial_{x_i}{V_3}~\nu_{j} \int_{\partial 
\Omega}a_2^{ij}(x)~\partial_{x_i}{V_2}~\nu_{j}\\ 
&-\int_{\partial 
D_2}a_2^{ij}(x)~\partial_{x_i}{V_3}~\nu_{j}\int_{\partial 
\Omega}a_2^{ij}(x)~\partial_{x_i}{V_1}~\nu_{j}, 
\end{split} 
\end{equation} 
then $Q_{\varepsilon}: C^2(\partial\Omega)\rightarrow \mathbb{R}$ is 
a linear functional. 
\begin{thm} \label{thm:lowbdd general} 
With the same conditions in Theorem \ref{thm:upbdd general}, let 
$u\in H^1(\Omega)\cap C^1(\overline{\widetilde{\Omega}})$ be the 
solution to equation (\ref{eq:k +infty general}). For $\varepsilon$ 
sufficiently small and $Q_{\varepsilon}[\varphi]$ defined by 
(\ref{eq:Q_e general}), there is a positive constant $C$ which 
depends only on $n$, $\kappa_0$, $r_0$, $\|\partial 
D_1\|_{C^{2,\alpha}}$, $\|\partial D_2\|_{C^{2,\alpha}}$, $\lambda$, 
$\Lambda$ and $\|A_2\|_{C^2(\overline{\widetilde\Omega})}$, but 
independent of $\varepsilon$ such that estimate (\ref{eq:lower}) 
holds. 
\end{thm} 
 
The paper is organized as follows. In Section 1 we prove Theorem 
\ref{thm:upbdd}$-$\ref{thm:lowbdd}. In Section 2 we give a criteria 
for $|Q_\varepsilon[\varphi]|$ to be bounded below by a positive 
constant independent of $\varepsilon$. Theorem \ref{thm:upbdd 
general}$-$\ref{thm:lowbdd general} are proved in Section 3. In the 
Appendix we present some elementary results for the conductivity 
problem.\\\\

%%%%%%%%%%%%%%%%%%%%%%%%%%%%%%%%%%%%%%%%%%%%%%%%%%%%%%%%%%%%%%%%%%%%%%%% 
% 
%      2. Proof of Theorem 
% 
%%%%%%%%%%%%%%%%%%%%%%%%%%%%%%%%%%%%%%%%%%%%%%%%%%%%%%%%%%%%%%%%%%%%%%%% 
 
\section {Proof of Theorem \ref{thm:upbdd} and \ref{thm:lowbdd}} 
 
In the introduction, we write $u=C_1v_1+C_2v_2+v_3$ as in 
(\ref{eq:decomp}). To prove our main theorems, we first estimate 
 $\|\nabla u\|_{L^{\infty}(\widetilde{\Omega})}$ 
in terms of  $|C_1-C_2|$, and then estimate $|C_1-C_2|$. \\\\ 
In this section we use, unless otherwise stated, $C$ to denote 
various positive constants whose values may change from line to line 
and which depend only on $n$, $\kappa_0$, $r_0$, $\|\partial 
\Omega\|_{C^{2,\alpha}}$, $\|\partial D_1\|_{C^{2,\alpha}}$ and 
$\|\partial D_2\|_{C^{2,\alpha}}$. 
\begin{prop} \label{prop:bd C_1-C_2} 
Under the hypotheses of Theorem \ref{thm:upbdd}, let $u$ be the 
solution of equation (\ref{eq:k +infty}). 
There exists a positive constants $C$, such 
that, for sufficiently small $\varepsilon>0$, 
\begin{equation} \label{ineq:bd:C_1-C_2} 
\frac{1}{\varepsilon}\mid C_{1}-C_{2} \mid\leq \|\nabla 
u\|_{L^{\infty}(\widetilde{\Omega})}\leq \frac{C}{\varepsilon}\mid 
C_{1}-C_{2} \mid + ~C\|\varphi\|_{C^2(\partial \Omega)}. 
\end{equation} 
\end{prop} 
%%%%%%%%%%%%%%%%%%%%%%%%%%%%%%%%%%%%%%%%%%%%%%%%%%%%%%%%%%%%%%%%%%%%%%%% 
% 
%      estimate of \nabla v_i 
% 
%%%%%%%%%%%%%%%%%%%%%%%%%%%%%%%%%%%%%%%%%%%%%%%%%%%%%%%%%%%%%%%%%%%%%%%% 
To prove this proposition, we first estimate the gradients of 
$v_1$, $v_2$ and $v_3$. Without loss of generality, we may assume 
throughout the proof of the proposition 
that $\|\varphi\|_{C^2(\partial \Omega)}=1$; see  Remark 
\ref{rem:varphi=1}. 
\begin{lem}\label{lm grad v_1,2} 
Let $v_1, v_2$ be defined by equations (\ref{eq:v1}) and 
({\ref{eq:v2}}), then for $n\geq 2$, we have 
\[\|\nabla 
v_1\|_{L^{\infty}(\widetilde{\Omega})}+\|\nabla 
v_2\|_{L^{\infty}(\widetilde{\Omega})}\leq 
\frac{C}{\varepsilon},~~~\|\frac{\partial v_1}{\partial 
\nu}\|_{L^{\infty}(\partial\Omega)}+\|\frac{\partial v_2}{\partial 
\nu}\|_{L^{\infty}(\partial\Omega)}\leq C.\] 
\end{lem} 
\emph{Proof}:\hspace{0.2cm} By the maximum principle, 
$\|v_1\|_{L^{\infty}(\widetilde{\Omega})}\leq 1$, and since $v_1$ 
achieves constants on each connected component of 
$\partial\widetilde \Omega$, and each connected component of 
$\partial\widetilde \Omega$ is $C^{2,\alpha}$ then the gradient 
estimates for harmonic functions implies that 
$$\|\nabla 
v_1\|_{L^{\infty}(\widetilde{\Omega})}\leq \frac 
{C\|v_1\|_{L^{\infty}}}{\text{dist}(\partial D_1,\partial 
D_2)}=\frac{C}{\varepsilon}.$$ Similarly, we can prove $\|\nabla 
v_2\|_{L^{\infty}(\widetilde{\Omega})}\leq C/\varepsilon$. The second 
inequality follows from the boundary estimates for harmonic functions.$\hfill\square$\\\\ 
Before estimating $|\nabla v_3|$, we first prove: 
\begin{lem} \label{lm grad rho} 
Let $\rho\in C^2(\widetilde{\Omega})$ be the solution to: 
\begin{equation} \label{eq: rho} 
        \left\{ \begin{aligned} 
           &\Delta \rho=0 \hspace{2.5cm}in~\widetilde{\Omega}, \\ 
           &\rho=0 ~~on~\partial D_1\cup\partial D_2,~~~~\rho=1 ~~on~\partial\Omega.\\ 
          \end{aligned} 
\right. 
\end{equation} 
Then $\|\nabla \rho\|_{L^{\infty}(\widetilde{\Omega})}\leq C$. 
\end{lem} 
\emph{Proof}:\hspace{0.2cm} Let $\rho_i(i=1,2)\in 
C^2(\Omega\backslash\overline{D}_i)\cap 
C^1(\overline{\Omega\backslash D_i})$ be the solution to: 
\begin{equation} 
        \left\{ \begin{aligned} 
           &\Delta \rho_i=0 \hspace{2.1cm}in~\Omega\backslash\overline{D}_i, \\ 
           &\rho_i=0 ~~on~\partial D_i,~~~~\rho_i=1 ~~on~\partial\Omega.\\ 
          \end{aligned} 
\right. \nonumber 
\end{equation} 
Again by the maximum principle and the strong maximum principle, we 
obtain $0<\rho_1<1$ in $\Omega\backslash\overline{D}_1$. Since 
$\overline{D}_2\subset\Omega\backslash\overline{D}_1$, we have 
$\rho_1>0=\rho$ on $\partial D_2$. And since $\rho_1=\rho$ on 
$\partial D_1$ and $\partial \Omega$, therefore $\rho_1>\rho$ on 
$\widetilde{\Omega}$. Now because $\rho_1=\rho=0$ on $\partial D_1$ 
and $\rho_1>\rho>0$ on $\widetilde{\Omega}$, so 
\[\|\nabla \rho\|_{L^{\infty}(\partial D_1)}\leq \|\nabla \rho_1\|_{L^{\infty}(\partial D_1)}\leq C.\] 
Similarly, 
\[\|\nabla \rho\|_{L^{\infty}(\partial D_2)}\leq \|\nabla \rho_2\|_{L^{\infty}(\partial D_2)}\leq C.\] 
By the boundary estimate of harmonic functions, we know that 
$\|\nabla \rho\|_{L^{\infty}(\partial \Omega)}\leq C$.\\ 
Since $\Delta\rho=0$ in $\widetilde{\Omega}$, $\partial_{x_i}\rho$ 
is also harmonic, by the maximum principle, 
\[\|\nabla 
\rho\|_{L^{\infty}(\widetilde{\Omega})}\leq\max\Big(\|\nabla 
\rho\|_{L^{\infty}(\partial D_1)},\|\nabla 
\rho\|_{L^{\infty}(\partial D_2)},\|\nabla 
\rho\|_{L^{\infty}(\partial \Omega)}\Big)\leq C. 
\]$\hfill\square$\\\\ 
Now, we estimate $|\nabla v_3|$: 
\begin{lem} \label{lm grad v3} 
Let $v_3$ be defined by equation (\ref{eq:v3}), for $n\geq 2$, we 
have 
\[\|\nabla v_3\|_{L^{\infty}(\widetilde{\Omega})}\leq C.\] 
\end{lem} 
\emph{Proof}:\hspace{0.2cm} Since $v_3=-\rho=\rho=0$ on $\partial 
D_i(i=1,2)$, and $-\rho\leq v_3=\varphi \leq \rho$ on 
$\partial\Omega$, we have, by the maximum principle, 
\[-\rho\leq v_3\leq \rho ~~\text{ in }\widetilde{\Omega}.\] 
It follows, for $i=1,2$, that 
\[\|\nabla 
v_3\|_{L^{\infty}(\partial D_i)}\leq \|\nabla 
\rho\|_{L^{\infty}(\partial D_i)}\leq C.\] 
By the boundary estimate, 
\[\|\nabla 
v_3\|_{L^{\infty}(\partial\Omega)}\leq C.\] By the harmonicity of 
$\partial_{x_i}v_3$ and the maximum principle, 
\[\|\nabla 
v_3\|_{L^{\infty}(\widetilde{\Omega})}\leq C.\]$\hfill\square$ 
 
\begin{rem} Without assuming 
$\|\varphi\|_{ C^2(\partial \Omega) }=1$, we have 
\[\|\nabla v_3\|_{L^{\infty}(\partial D_1\cup \partial D_2)}\leq C 
\|\varphi\|_{L^{\infty}(\partial\Omega)},\] where $C$ has the 
dependence specified at the beginning of this section, except that it 
does not depend on $\|\partial \Omega\|_{ C^{2, \alpha} }$. 
This is 
easy to see from the proof of Lemma \ref{lm grad v3}. 
 \label{rem1.1} 
\end{rem} 
 The above lemma yields 
the main result of \cite{ADKL}. 
\begin{cor} \label{cor grad v3} 
\emph{(\cite{ADKL})} Let $B_1$ and $B_2$ be two spheres with radius 
$R$ and centered at $(\pm R\pm \frac{\varepsilon}{2},0,\cdots,0)$, 
respectively. Let $H$ be a harmonic function in $\mathbb{R}^3$. 
Define $u$ to be the solution to 
\begin{equation} 
        \left\{ \begin{aligned} 
           &\Delta u=0 ~~~~~~~~~~~~~~~~~~~~~~~~~in~\mathbb{R}^3\backslash\overline{B_1\cup B_2}, \\ 
           &u=0 ~~~~~~~~~~~~~~~~~~~~~~~~~~~on~\partial B_1\cup\partial B_2, \\ 
           &u(x)-H(x)=O(|x|^{-1}) ~~~~~as~|x|\rightarrow +\infty. 
          \end{aligned} 
\right. 
\nonumber 
\end{equation} 
Then there is a constant $C$ independent of $\varepsilon$ such that 
$$\|\nabla (u-H)\|_{L^{\infty}(\mathbb{R}^3\backslash\overline{B_1\cup B_2})}\leq C.$$ 
\end{cor} 
\emph{Proof}:\hspace{0.2cm} By the maximum principle and interior 
estimates of harmonic functions, the $C^3$ norm of  $u|_{B_{2R}(0)}$ 
is bounded by a constant independent of $\varepsilon$. Apply Lemma 
\ref{lm grad v3} with $\Omega= B_{2R}(0)$ and 
$\varphi=u|_{B_{2R}(0)}$, we immediately obtain 
the above corollary.$\hfill\square$\\ 
% 
% 
%\nabla v_1, v_2, v_3 
\\ 
With the above lemmas, we give the\\\\ 
\emph{Proof of Proposition \ref{prop:bd C_1-C_2}}:\hspace{0.2cm} 
Since $u=C_1$ on $\partial D_1$, $u=C_2$ on $\partial D_2$, 
$\text{dist}(\partial D_1,\partial D_2)=\varepsilon$, by the mean 
value theorem, $\exists ~\xi\in \widetilde{\Omega}$ such that 
\[\|\nabla 
u\|_{L^{\infty}(\widetilde{\Omega})}\geq|\nabla u(\xi)|\geq 
\frac{|C_1-C_2|}{\varepsilon}.\] By the 
decomposition formula 
(\ref{eq:decomp}), 
\[\nabla u=C_1\nabla v_1+C_2\nabla v_2+\nabla v_3=(C_1-C_2)\nabla v_1+C_2\nabla(v_1+v_2)+\nabla v_3.\] 
Hence, 
\[\|\nabla 
u\|_{L^{\infty}(\widetilde{\Omega})}\leq |C_1-C_2|\|\nabla 
v_1\|_{L^{\infty}(\widetilde{\Omega})}+|C_2|\|\nabla 
(v_1+v_2)\|_{L^{\infty}(\widetilde{\Omega})}+\|\nabla 
v_3\|_{L^{\infty}(\widetilde{\Omega})}.\] By Lemma \ref{lm grad 
rho}, since $v_1+v_2=1-\rho$ in $\widetilde{\Omega}$, we have 
\[\|\nabla(v_1+v_2)\|_{L^{\infty}(\widetilde{\Omega})}=\|\nabla(1-\rho)\|_{L^{\infty}(\widetilde{\Omega})}=\|\nabla\rho\|_{L^{\infty}(\widetilde{\Omega})}\leq C.\] 
Using the fact we showed in the Appendix, $\|u\|_{H^1(\Omega)}\leq 
C$, so $|C_1|+|C_2|\leq C$.\\ 
Therefore using also Lemma \ref{lm grad v_1,2} we obtain, 
\[\|\nabla 
u\|_{L^{\infty}(\widetilde{\Omega})}\leq \frac{C}{\varepsilon}\mid 
C_1-C_2 \mid + ~C.\]This proof is 
now completed.$\hfill\square$\\\\ 
Later we will give an estimate of $|C_1-C_2|$, which, 
 together with 
Proposition \ref{prop:bd C_1-C_2}, 
 yields the lower and upper bounds 
of $\|\nabla u\|_{L^\infty(\widetilde\Omega)}$ for strictly convex 
subdomains $D_1$ 
and $D_2$.\\\\ 
%%%%%%%%%%%%%%%%%%%%%%%%%%%%%%%%%%%%%%%%%%%%%%%%%%%%%%%%%%%%%%%%%%%%%%%% 
% 
%      Definition and Propositions of a_ii and b_i 
% 
%%%%%%%%%%%%%%%%%%%%%%%%%%%%%%%%%%%%%%%%%%%%%%%%%%%%%%%%%%%%%%%%%%%%%%%% 
\subsection{Estimate of $|C_{1}-C_{2}|$} 
Back to the decomposition formula (\ref{eq:decomp}), denote 
\begin{equation} \label{eq:ab} 
a_{ij}=\int_{\partial D_i}\frac{\partial 
v_j}{\partial\nu}\hspace{.5cm}(i,j=1,2), \hspace{.5cm} 
b_i=\int_{\partial D_i}\frac{\partial 
v_3}{\partial\nu}\hspace{.5cm}(i=1,2). 
\end{equation} 
We first give some basic lemmas: 
\begin{lem}\label{lm:ab} 
Let $a_{ij}$ and $b_{i}$ be defined as in (\ref{eq:ab}), then they 
satisfy the following: 
\begin{enumerate} 
   \item   $a_{12}=a_{21}>0,\  a_{11}<0,\  a_{22}<0$, 
   \item   $-C\leq a_{11}+a_{21}\leq -\frac{1}{C}, \ -C\leq a_{22}+a_{12}\leq -\frac{1}{C}$, 
   \item   $|b_1|\leq C, ~|b_2|\leq C$. 
\end{enumerate} 
\end{lem} 
By the fourth line of equation (\ref{eq:k +infty}), $C_1$ and $C_2$ 
satisfy 
\begin{equation} \label{eq:c1c2} 
\left\{ \begin{aligned} 
           &a_{11} C_1 + a_{12} C_2 + b_1=0, \\ 
           &a_{21} C_1 + a_{22} C_2 + b_2=0. 
          \end{aligned} 
\right. 
%\nonumber 
\end{equation} 
By solving the above linear system, using $a_{12}=a_{21}$ and 
$a_{11}a_{22}-a_{12}a_{21}>0$ which follows from Lemma \ref{lm:ab}, 
we obtain 
\begin{align}\label{eq:c1c2s} 
C_1=\frac{-b_1a_{22}+b_2a_{12}}{a_{11}a_{22}-a_{12}^2}, \ \ \ \ 
C_2=\frac{-b_2a_{11}+b_1a_{12}}{a_{11}a_{22}-a_{12}^2}, 
\end{align} 
and therefore, 
\begin{equation}\label{eq:c1-c2} 
|C_1-C_2|=\frac{|b_1-\alpha b_2|}{|a_{11}-\alpha a_{12}|} , 
\hspace{1cm} \text{where} \ 
~\alpha=\frac{a_{11}+a_{12}}{a_{22}+a_{12}}>0. 
\end{equation} 
Based on this formula, we will give the estimates for 
$|a_{11}-\alpha a_{12}|$ and $|b_1-\alpha b_2|$, then the estimate 
for $|C_1-C_2|$ follows immediately.\\\\ 
\emph{Proof of Lemma \ref{lm:ab}}:\hspace{0.2cm} (1) By the maximum 
principle and the strong maximum principle, 
$$0<v_1<1 ~~~\text{in} ~\widetilde\Omega.$$ 
By the Hopf Lemma, we know that 
$$\frac{\partial 
v_1}{\partial\nu}\big|_{\partial D_1}<0, ~~~\frac{\partial 
v_1}{\partial\nu}\big|_{\partial D_2}>0, ~~~\frac{\partial 
v_1}{\partial\nu}\big|_{\partial \Omega}<0.$$ 
Similarly, 
$$\frac{\partial v_2}{\partial\nu}\big|_{\partial D_1}>0, 
~~~\frac{\partial v_2}{\partial\nu}\big|_{\partial D_2}<0, 
~~~\frac{\partial v_2}{\partial\nu}\big|_{\partial \Omega}<0.$$ Thus 
$a_{11}<0$, $a_{12}>0$, $a_{21}>0$ and $a_{22}<0$. 
 
Also, since $v_1$ and $v_2$ are the solutions of equations 
(\ref{eq:v1}) and equations (\ref{eq:v2}), respectively, we have 
\begin{equation} \label{eq:a12=a21} 
\begin{split} 
0&=\int_{\widetilde\Omega}\Delta v_1\cdot v_2 - 
\int_{\widetilde\Omega}\Delta v_2\cdot v_1 = -\int_{\partial 
D_2}\frac{\partial v_1}{\partial\nu} \cdot 1 + \int_{\partial 
D_1}\frac{\partial v_2}{\partial\nu} \cdot 1\\ 
&=-a_{21}+a_{12}, 
\end{split} 
\end{equation} 
i.e. $a_{21}=a_{12}$.\\ 
 
(2) We will prove the first inequality, the second one stands with 
the same reason. By the harmonicity of $v_1$ in $\widetilde\Omega$, 
$$a_{11}+a_{21}=-\int_{\widetilde\Omega}\Delta v_1+\int_{\partial \Omega}\frac{\partial 
v_1}{\partial\nu}=\int_{\partial \Omega}\frac{\partial 
v_1}{\partial\nu}<0.$$  By Lemma \ref{lm grad v_1,2}, 
$$a_{11}+a_{21}=\int_{\partial \Omega}\frac{\partial 
v_1}{\partial\nu}\geq -C.$$ 
On the other hand, since $0<v_1<1 
~in~\widetilde\Omega$ and $v_1=1 ~ on~ \partial D_1$, by the 
boundary gradient estimates of a harmonic function, $\exists ~B(\bar 
x, 2\bar r)\subset\widetilde\Omega$, such that $v_1>1/2$ in $B(\bar 
x, \bar r)$, where $\bar r$ is independent of $\varepsilon$. Let 
$\rho\in C^2(\Omega\backslash\overline {D_2\cup B(\bar x, \bar r)}) 
\cup C^1(\partial \Omega \cup \partial D_2 \cup 
\partial B(\bar x, \bar r))$ 
be the solution of the following equation: 
\begin{equation} 
        \left\{ \begin{aligned} 
           &\Delta \rho=0 \hspace{2.5cm}in~\Omega\backslash\overline {D_2\cup B(\bar x, \bar r)}, \\ 
           &\rho=1/2 ~~on~\partial B(\bar x, \bar r)~~~\rho=0 ~~on~\partial D_2\cup\partial \Omega.\\ 
          \end{aligned} 
\right. \nonumber 
\end{equation} 
By the maximum principle and the strong maximum principle, 
$0<\rho<1/2$ in $\Omega\backslash\overline{D_2\cup B(\bar x, \bar 
r)}.$ A contradiction argument based on the Hopf Lemma yields, 
$$-\frac{\partial 
\rho}{\partial\nu}\geq \frac{1}{C} ~~~~~\text{on}~ 
\partial\Omega.$$ 
On the other hand, since $\rho\leq v_1$ on the boundary of 
$\Omega\backslash\overline{D_1\cup D_2\cup B(\bar x, \bar r)},$ 
we obtain, via 
the maximum principle, $0<\rho\leq v_1$ in 
$\Omega\backslash\overline{D_1\cup D_2\cup B(\bar x, \bar r)}.$ It 
follows, using $\rho=v_1=0$ on $\partial\Omega$, that 
$$\frac{\partial 
v_1}{\partial\nu}\leq\frac{\partial 
\rho}{\partial\nu}~~~~~\text{on}~ 
\partial\Omega.$$ 
Thus, 
$$a_{11}+a_{21}=\int_{\partial 
\Omega}\frac{\partial v_1}{\partial\nu}\leq\int_{\partial 
\Omega}\frac{\partial \rho}{\partial\nu}\leq-\frac{1}{C}.$$ 
 
(3) Clearly, 
\begin{equation} 
0=\int_{\widetilde\Omega}\Delta v_1\cdot v_3 - 
\int_{\widetilde\Omega}\Delta v_3\cdot v_1 = \int_{\partial 
\Omega}\frac{\partial v_1}{\partial\nu} \cdot \varphi + 
\int_{\partial D_1}\frac{\partial v_3}{\partial\nu} \cdot 
1=\int_{\partial \Omega}\frac{\partial v_1}{\partial\nu} \cdot 
\varphi+b_1. \nonumber 
\end{equation} 
Thus, 
$$|b_1|=\Big|\int_{\partial \Omega}\frac{\partial v_1}{\partial\nu} 
\cdot \varphi\Big|\leq \int_{\partial \Omega}\Big|\frac{\partial 
v_1}{\partial\nu}\Big| \leq C.$$Thus, we finished the 
proof.$\hfill\square$\\\\

\subsection{Estimate of $|a_{11}-\alpha a_{12}|$} 
By a translation and rotation of the axis, we may assume without 
loss of generality that $D_1$, $D_2$ are two strictly convex 
subdomains in $\Omega\subset \mathbb{R}^n$ which satisfy the 
following: 
\begin{equation}\label{D1, D2} 
(-\varepsilon/2, 0')\in 
\partial D_1, ~(\varepsilon/2, 0')\in \partial D_2,~\varepsilon=\text{dist}(\partial D_1, \partial D_2)=\text{dist}(D_1, 
D_2). 
\end{equation} 
%When $\varepsilon=0$, two inclusions touch each other at the origin, 
%and we denote that $D_1$ becomes $D_1^*$ and $D_2$ becomes $D_2^*$. 
%Therefore, we can also consider that $D_1$ is obtained by moving 
%$D^*_1$ to the left for $\varepsilon/2$ and $D_2$ is obtained by 
%moving $D^*_2$ to the right for $\varepsilon/2$. Now, we also assume 
%$\partial D^*_1$, $\partial D^*_2$ are $C^{2,\alpha}$ and 
Near the origin, we can find a ball $B(0,r)$ such that the portion 
of $\partial D_i~(i=1,2)$ in $B(0,r)$ is strictly convex, where 
$r>0$ is independent of $\varepsilon$. Then $\partial D_1\cap 
B(0,r)$ and $\partial D_2\cap B(0,r)$ can be represented by the 
graph of $x_1=f(x')-\varepsilon/2$ and $x_1=g(x')+\varepsilon/2$ 
respectively, where $x'=(x_2,\cdots, x_n)$. Thus $f(0')=g(0')=0, 
~\nabla f(0')=\nabla g(0')=0$, and $-CI\leq\big(D^2f(0')\big)\leq 
-\frac{1}{C}I$, $\frac{1}{C}I\leq\big(D^2g(0')\big)\leq CI$.\\\\ 
With these notations, we first estimate $a_{ii}$ for $i=1,2$. 
\begin{lem}\label{lm: bd aii n 2} 
Let $a_{ii}$ be defined by (\ref{eq:ab}), then 
$$\frac{1}{C\sqrt{\varepsilon}}\leq-a_{ii}\leq\frac{C}{\sqrt{\varepsilon}}, ~~~for ~n=2, ~i=1,2.$$ 
\end{lem} 
\emph{Proof}:\hspace{0.2cm} It suffices to prove it for $a_{11}$. By 
the harmonicity of $v_1$, we have 
\begin{equation} 
0=\int_{\widetilde\Omega}\Delta v_1 \cdot 
v_1=-\int_{\widetilde\Omega}|\nabla v_1|^2 -\int_{\partial 
D_1}\frac{\partial v_1}{\partial\nu}=-\int_{\widetilde\Omega}|\nabla 
v_1|^2-a_{11}, 
\nonumber 
\end{equation} 
i.e. 
$$a_{11}=-\int_{\widetilde\Omega}|\nabla 
v_1|^2.$$ 
 
Now we construct a function (here in $\mathbb{R}^2$, we let $x=x_1, 
~y=x_2$) 
\begin{equation}\label{eq:w n=2} 
\overline{w}(x,y)=-\frac{x-g(y)-\frac{\varepsilon}{2}}{g(y)-f(y)+\varepsilon} 
\end{equation} 
on $O_{r}:=\widetilde\Omega\cap{\left\{(x,y)\big |~|y|<r\right\}}$. 
It is clear that $\overline{w}(x,y)$ is linear in x for fixed y and 
$$\overline{w}\mid_{B(0,r)\cap\partial 
D_1}=1;~~~\overline{w}\mid_{B(0,r)\cap\partial D_2}=0,$$ 
so we have 
$$\int^{g(y)+\frac{\varepsilon}{2}}_{f(y)-\frac{\varepsilon}{2}}|\partial_{x}\overline{w}(x,y)|^2dx\leq 
\int^{g(y)+\frac{\varepsilon}{2}}_{f(y)-\frac{\varepsilon}{2}}|\partial_{x}v_1(x,y)|^2dx,$$ 
 
i.e. 
$$\frac{1}{g(y)-f(y)+\varepsilon}\leq\int^{g(y)+\frac{\varepsilon}{2}}_{f(y)-\frac{\varepsilon}{2}}|\partial_{x}v_1(x,y)|^2.$$ 
 
Integrating on y we get 
\begin{equation}\label{eq:v1xlow} 
\begin{split} 
&\int^{r/2}_{0}\int^{g(y)+\frac{\varepsilon}{2}}_{f(y)-\frac{\varepsilon}{2}}|\partial_{x}v_1(x,y)|^2dxdy 
\geq\int^{r/2}_{0}\frac{1}{g(y)-f(y)+\varepsilon}dy\\ 
\geq&~\frac 
1C\int^{r/2}_{0}\frac{1}{y^2+\varepsilon}dy=\frac{1}{C\sqrt\varepsilon}. 
\end{split} 
\end{equation} 
Thus 
$$-a_{11}\geq\int^{r/2}_{0}\int^{g(y)+\frac{\varepsilon}{2}}_{f(y)-\frac{\varepsilon}{2}}|\partial_{x}v_1(x,y)|^2dxdy\geq\frac{1}{C\sqrt\varepsilon}.$$ 
On the other hand, we can find $\psi\in C^2(\overline{\Omega})$ such 
that 
$$\psi=0 ~on~\overline O_{r/8}, ~~\psi=1 ~on ~\partial 
D_1\backslash(\overline{O_{r/4}}), ~~\psi=0 ~on~ \partial 
D_2\backslash(\overline{O_{r/4}}), 
$$ 
$$\psi=0 ~on ~\partial\Omega,~~~\text{and}~~~\|\nabla \psi\|_{L^\infty(\Omega)} 
\le C.$$ 
We can also find $\rho\in C^2(\overline{\Omega})$ such that 
$$0\leq\rho\leq1,~ \rho=1 ~on~ \overline O_{r/2},~ \rho=0 ~on~ 
\overline\Omega\backslash O_{r} ~\text{and}~ |\nabla\rho|\leq C.$$ 
Let $w=\rho\overline{w}+(1-\rho)\psi$, then $w=1=v_1$ on $\partial 
D_1$;$w=0=v_1$ on $\partial D_2$; $w=0=v_1$ on $\partial \Omega$ and 
$w=\overline{w}$ on $\overline O_{r/2}$. Then by the properties of 
$\psi$, $\rho$ and the harmonicity of $v_1$, we have 
\begin{equation} 
\label{eq:v1w} \int_{\widetilde\Omega}|\nabla 
v_1|^2\leq\int_{\widetilde\Omega}|\nabla w|^2 
\leq\int_{\widetilde\Omega\cap O_{r/2}}|\nabla \overline w|^2+C. 
\end{equation} 
A calculation gives 
$$\partial_y 
\overline{w}=\frac{g'(y)(g(y)-f(y)+\varepsilon)-(g(y)-x+\frac{\varepsilon} 
{2})(g'(y)-f'(y))}{(g(y)-f(y)+\varepsilon)^2}.$$ 
We will show $\int_{\widetilde \Omega \cap O_{r/2}}|\partial_y \overline{w}|^2\leq C$. \\ 
Indeed, 
\begin{equation}\label{eq:vwy} 
\begin{split} 
&\int^{r/2}_{0}\int^{g(y)+\frac{\varepsilon}{2}}_{f(y)-\frac{\varepsilon}{2}}|\partial_{y}\overline{w}(x,y)|^2dxdy\\ 
&\leq 
2\int^{r/2}_{0}\int^{g(y)+\frac{\varepsilon}{2}}_{f(y)-\frac{\varepsilon}{2}} 
\left(\frac{g'(y)^2}{(g(y)-f(y)+\varepsilon)^2}+ 
\frac{(g(y)-x+\frac{\varepsilon} 
{2})^2(g'(y)-f'(y))^2}{(g(y)-f(y)+\varepsilon)^4}\right)dxdy\\ 
&=2\int^{r/2}_{0}\frac{g'(y)^2}{g(y)-f(y)+\varepsilon}dy+2\int^{r/2}_{0}\frac{(g'(y)-f'(y))^2}{g(y)-f(y)+\varepsilon}dy\\ 
&\leq C\int^{r/2}_{0}\frac{y^2}{y^2+\varepsilon}dy+C\int^{r/2}_{0}\frac{y^2}{y^2+\varepsilon}dy\\ 
&\leq C. 
\end{split} 
\end{equation} 
Then by (\ref{eq:v1w}) and (\ref{eq:vwy}) 
\begin{equation} 
\begin{split}\label{eq:v1xup} 
|a_{11}|&=\int_{\widetilde\Omega}|\nabla v_1|^2\leq\int_{\widetilde\Omega\cap O_{r/2}}|\nabla \overline w|^2+C\\ 
&\leq C\int^{r/2}_{0}\int^{g(y)+\frac{\varepsilon}{2}}_{f(y)-\frac{\varepsilon}{2}}|D_{x}\overline{w}(x,y)|^2dxdy+C\\ 
&=C\int^{r/2}_{0}\frac{1}{g(y)-f(y)+\varepsilon}dy+C\leq C\int^{r/2}_{0}\frac{1}{y^2+\varepsilon}dy+C\\ 
&\leq \frac{C}{\sqrt\varepsilon}. 
\end{split} 
\end{equation} 
The proof is completed.$\hfill\square$\\\\ 
Similarly, we have 
\begin{lem}\label{lm: bd aii n 3} 
Let $a_{ii}$ be defined by (\ref{eq:ab}), 
$$\frac{1}{C}|\ln 
\varepsilon|\leq-a_{ii}\leq C|\ln\varepsilon|, ~~ for ~ n=3, 
~i=1,2.$$ 
\end{lem} 
\emph{Proof}:\hspace{0.2cm} We consider 
\begin{equation}\label{eq:w n>2} 
\overline{w}(x_1,x')=-\frac{x-g(x')-\frac{\varepsilon}{2}}{g(x')-f(x')+\varepsilon} 
\end{equation} 
on $O_{r/2}:=\widetilde\Omega\cap\{(x_1,x')|~|x'|<\frac{r}{2}\}$. 
Use the same proof in Lemma \ref{lm: bd aii n 2}, we have 
$$\int^{r/2}_{0}\int^{g(x')+\frac{\varepsilon}{2}}_{f(x')-\frac{\varepsilon}{2}}|\partial_{x'}\overline{w}(x_1,x')|^2dx_1dx'\leq C.$$ 
Therefore, it suffices to verify that 
$$\int_{\widetilde \Omega \cap O_{r/2}}|\partial_{x_1}\overline{w}(x_1,x')|^2\sim |\ln \varepsilon|.$$ 
Indeed, 
$$\int_{\widetilde \Omega \cap O_{r/2}}|\partial_{x_1}\overline{w}(x_1,x')|^2=\int_{|x'|<r/2}\frac{1}{g(x')-f(x')+\varepsilon}dx'\sim\int_{0}^{r/2}\frac{t}{Ct^2+\varepsilon}dt\sim|\ln\varepsilon|.$$ 
This completes the proof.$\hfill\square$\\ 
 
\begin{lem}\label{lm: bd aii n 4+} 
Let $a_{ii}$ be defined by (\ref{eq:ab}), 
$$\frac{1}{C}\leq-a_{ii}\leq C ~~~ for~ n\geq 4, ~i=1,2.$$ 
\end{lem} 
\emph{Proof}:\hspace{0.2cm} We only need 
$$\int_{O_{r/2}}|\partial_{x_1}\overline{w}(x_1,x')|^2=\int_{|x'|<r/2}\frac{1}{g(x')-f(x')+\varepsilon}dx'\sim\int_{0}^{r/2}\frac{t^{n-2}}{Ct^2+\varepsilon}dt\sim C.$$ 
The proof is completed.$\hfill\square$\\ 
 
\begin{lem}\label{lm:bd alpha} 
Let $\alpha$ be defined by (\ref{eq:c1-c2}), we have 
$$\frac{1}{C}\leq \alpha\leq C.$$ 
\end{lem} 
\emph{Proof}:\hspace{0.2cm} By the definition of $\alpha$ and using 
the second statement in Lemma \ref{lm:ab}, we are 
done.$\hfill\square$\\\\ 
To summarize, we have 
\begin{prop}\label{prop: bd a11-alpha a22} 
Let $a_{ij}$ and $\alpha$ be defined by (\ref{eq:ab}) and 
(\ref{eq:c1-c2}), we have 
\begin{enumerate} 
   \item   $\frac{1}{C\sqrt\varepsilon}\leq|a_{11}-\alpha 
   a_{12}|\leq\frac{C}{\sqrt\varepsilon} ~~~~~for~ n=2,$ 
   \item  $\frac{1}{C}|\ln\varepsilon|\leq|a_{11}-\alpha 
   a_{12}|\leq C|\ln\varepsilon| ~~~~~for~ n=3,$ 
   \item   $\frac{1}{C}\leq|a_{11}-\alpha 
   a_{12}|\leq C ~~~~~for~ n\geq 4.$ 
\end{enumerate} 
\end{prop} 
\emph{Proof}:\hspace{0.2cm} Since $a_{11}<0$, $a_{12}>0$, 
$a_{11}+a_{12}<0$ and $\alpha>0$, we have 
$$|a_{11}|<|a_{11}-\alpha 
a_{12}|<(1+\alpha)|a_{11}|.$$ Combining the results of Lemma 
\ref{lm: bd aii n 2}, Lemma \ref{lm: bd aii n 3}, Lemma \ref{lm: bd 
aii n 4+} and Lemma \ref{lm:bd alpha}, the proof is 
completed.$\hfill\square$\\ 
 
\subsection{Estimate of $|b_{1}-\alpha b_{2}|$} 
\begin{prop}\label{prop:bd b1-alpha b2} 
Let $b_1, ~b_2$, $\alpha$ and $Q_\varepsilon[\varphi]$ be defined by 
(\ref{eq:ab}), (\ref{eq:c1-c2}) and (\ref{eq:Q_e}), we have 
$$\frac{|Q_\varepsilon[\varphi]|}{C}\leq|b_{1}-\alpha b_{2}|\leq C\|\varphi\|_{C^2(\partial \Omega)}.$$ 
\end{prop} 
\emph{Proof}:\hspace{0.2cm} Combining the third result in Lemma 
\ref{lm:ab} and Lemma \ref{lm:bd alpha}, we have 
$$|b_{1}-\alpha b_{2}|\leq |b_1|+|\alpha||b_2|\leq C\|\varphi\|_{C^2(\partial \Omega)}.$$ 
On the other hand, by the definition and the harmonicity of $v_1$ 
and $v_2$ and using Lemma \ref{lm:ab}, we obtain 
\begin{equation} 
\begin{split} 
|b_{1}-\alpha b_{2}|&=\frac{|b_1(a_{22}+a_{12})-b_2(a_{11}+a_{12})|}{|a_{22}+a_{12}|}\\ 
&\geq \frac{1}{C}\cdot\Big|\int_{\partial D_1}\frac{\partial 
v_3}{\partial \nu}\int_{\partial \Omega}\frac{\partial v_2}{\partial 
\nu}-\int_{\partial D_2}\frac{\partial v_3}{\partial 
\nu}\int_{\partial \Omega}\frac{\partial v_1}{\partial 
\nu}\Big|=\frac{|Q_\varepsilon[\varphi]|}{C}.\\ 
\end{split} 
\nonumber 
\end{equation} 
This completes the proof.$\hfill\square$\\\\ 
Now we are ready to prove our two main theorems:\\\\ 
\emph{Proof of Theorem 
\ref{thm:upbdd}-\ref{thm:lowbdd}}:\hspace{0.2cm} By Proposition 
\ref{prop:bd C_1-C_2} and (\ref{eq:c1-c2}), then using Proposition 
\ref{prop: bd a11-alpha a22}, \ref{prop:bd b1-alpha b2}, we are 
done.$\hfill\square$ 

\bigskip

As we mentioned in Remark \ref{rem1}, the strict convexity 
assumption of the two inclusions can be weakened. In fact, our 
proofs of Theorem \ref{thm:upbdd}$-$\ref{thm:lowbdd} apply, with 
minor modification, to more general situations: 
 
In $\mathbb{R}^n$, $n\geq 2$, under the same assumptions in the 
beginning of Section 1.2 except for the strict convexity 
condition, $\partial D_1\cap B(0,r)$ and $\partial D_2\cap B(0,r)$ 
can be represented by the graph of $x_1=f(x')-\frac{\varepsilon}{2}$ 
and $x_1=g(x')+\frac{\varepsilon}{2}$, then $f(0')=g(0')=0$, $\nabla 
(g-f)(0')=0$. Assume further that 
\begin{equation} \label{expan:g-f} 
\lambda_0|x'|^{2m}\leq g(x')-f(x')\leq \lambda_1|x'|^{2m}, ~~~ \forall |x'|\leq r/2,
\end{equation} 
for some
$\varepsilon$-independent $\lambda_0, \lambda_1>0, m\geq 1 \in \mathbb{Z}$. 
 
Under the above assumption, let $u\in H^1(\Omega)\cap 
C^1(\overline{\widetilde{\Omega}})$ be the solution to equation 
(\ref{eq:k +infty}). For $\varepsilon$ sufficiently small, there 
exist positive constants $C$ and $C'$, such that 
\begin{equation} \label{bd:R^n} 
\begin{split} 
\frac{|Q_{\varepsilon}[\varphi]|}{C'}\cdot\varepsilon^{-\frac{n-1}{2m}}\leq 
&\|\nabla u\|_{L^{\infty}(\widetilde{\Omega})}\leq 
C\|\varphi\|_{C^2(\partial 
\Omega)}\cdot\varepsilon^{-\frac{n-1}{2m}}, 
~~~~\emph{if}~ n-1<2m,\\ 
\frac{|Q_{\varepsilon}[\varphi]|}{C'}\cdot\frac{1}{\varepsilon|\ln\varepsilon|}\leq 
&\|\nabla u\|_{L^{\infty}(\widetilde{\Omega})}\leq 
C\|\varphi\|_{C^2(\partial 
\Omega)}\cdot\frac{1}{\varepsilon|\ln\varepsilon|}, 
~~~\emph{if}~ n-1=2m,\\ 
\frac{|Q_{\varepsilon}[\varphi]|}{C'}\cdot\frac{1}{\varepsilon}\leq 
&\|\nabla u\|_{L^{\infty}(\widetilde{\Omega})}\leq 
C\|\varphi\|_{C^2(\partial \Omega)}\cdot\frac{1}{\varepsilon}, 
~~~~~~~~~\emph{if}~ n-1>2m, 
\end{split} 
\end{equation} 
where $Q_{\varepsilon}[\varphi]$ is defined by (\ref{eq:Q_e}), and 
$C$ depends on $n$, $m$, $\lambda_0$, $\lambda_1$, $r_0$, $\|\partial 
\Omega\|_{C^{2,\alpha}}$, $\|\partial D_1\|_{C^{2,\alpha}}$ and 
$\|\partial D_2\|_{C^{2,\alpha}}$, $C'$ depends on the same as $C$ 
and also $\|\varphi\|_{C^2(\partial \Omega)}$, but both are 
independent of $\varepsilon$. 
 
The proof is essentially the same except for the computation of 
$\int_{\widetilde{\Omega}}|\nabla v_1|^2$.\\ 
In fact, 
$$\int^{r/2}_{0}\int^{g(x')+\frac{\varepsilon}{2}}_{f(x')-\frac{\varepsilon}{2}}|\partial_{x'}\overline{w}(x_1,x')|^2dx_1dx'\leq C,$$ 
still holds. Then by (\ref{eq:v1xlow}) and (\ref{eq:v1xup}) we only 
need to calculate 
$$\int_{|x'|<r/2}\frac{1}{g(x')-f(x')+\varepsilon}dx'\sim\int^{r/2}_{0}\frac{\rho^{n-2}}{\rho^{2m}+\varepsilon}d\rho.$$ 
Indeed, if $~n-1<2m$, 
$$\int^{r/2}_{0}\frac{\rho^{n-2}}{\rho^{2m}+\varepsilon}d\rho= 
\varepsilon^{\frac{n-1}{2m}-1}\int^{r/{2\varepsilon^{\frac{1}{2m}}}}_{0}\frac{s^{n-2}}{s^{2m}+1}ds\sim 
C\varepsilon^{\frac{n-1}{2m}-1},$$ if $~n-1=2m$, 
$$\int^{r/2}_{0}\frac{\rho^{n-2}}{\rho^{2m}+\varepsilon}d\rho= 
\frac{1}{2m}\int^{r/2}_{0}\frac{1}{\rho^{2m}+\varepsilon}d\rho^{2m}\sim 
C|\ln \varepsilon|,$$ if $~n-1>2m$, 
$$\int^{r/2}_{0}\frac{\rho^{n-2}}{\rho^{2m}+\varepsilon}d\rho\sim 
C.$$ Therefore, we obtain (\ref{bd:R^n}) by using the same arguments 
in the proofs of Theorem \ref{thm:upbdd} and Theorem 
\ref{thm:lowbdd}. 
 
Actually, we can replace $2m$ by any real number
$\beta>0$, the results still hold.

\section{Estimate of $|Q_\varepsilon[\varphi]|$} 
In order to identify situations when $\|\nabla u\|_{L^\infty}$ 
behaves exactly as the upper bound established in Theorem 
\ref{thm:upbdd}, 
 we  estimate 
 in this section 
$|Q_\varepsilon[\varphi]|$. 
To emphasize the dependence 
 on $\varepsilon$, we denote $D_1$, $D_2$ by 
$D_{1\varepsilon}$, $D_{2\varepsilon}$, 
denote $\varphi$ by $\varphi_\varepsilon$, 
and  denote $v_{1}$, $v_{2}$, 
$v_{3}$ defined by equation (\ref{eq:v1}), (\ref{eq:v2}), 
(\ref{eq:v3}) as $v_{1\varepsilon}$, $v_{2\varepsilon}$, 
$v_{3\varepsilon}$. 
In this section we assume, in addition to the 
hypotheses in Theorem \ref{thm:upbdd}, 
 that along a sequence 
$\varepsilon\to 0$ (we still denote it 
as $\varepsilon$), 
 $D_{1\varepsilon}\to D^*_1$, $D_{2\varepsilon}\to 
D^*_2$ in $C^{2, \alpha}$ norm, $\varphi_\varepsilon\to \varphi^*$ 
in $C^{1, \alpha}(\partial \Omega)$. 
We use notation 
$\widetilde\Omega^*=\Omega\backslash\overline{D_1^*\cup D_2^*}$, and 
 assume, without loss of generality, that $D^*_1\cap D^*_2=\{0\}$. 
We will show that as $\varepsilon\to 0$, 
$v_{i\varepsilon}$ converges, in appropriate sense, 
to  $v^*_i$ which satisfies 
\begin{equation} 
\label{eq:V1}\left\{\begin{aligned} 
             &\Delta v_1^*=0 ~~~~~~~~~~~~in~\widetilde\Omega^*, \\ 
             &v_1^*=1~~on~\partial D_1^*\backslash\{0\},~~~~v_1^*=0~~on~\partial\Omega\cup\partial D_2^*\backslash\{0\},\\ 
\end{aligned} 
\right. 
\end{equation} 
\begin{equation} 
\label{eq:V2}\left\{\begin{aligned} 
             &\Delta v_2^*=0 ~~~~~~~~~~~~in~\widetilde\Omega^*, \\ 
             &v_2^*=1~~on~\partial D_2^*\backslash\{0\},~~~~v_2^*=0~~on~\partial\Omega\cup\partial D_1^*\backslash\{0\},\\ 
\end{aligned} 
\right. 
\end{equation} 
 
\begin{equation} 
\label{eq:V3}\left\{\begin{aligned} 
             &\Delta v_3^*=0 ~~~~~~~~~~~~in~\widetilde\Omega^*,\\ 
             &v_3^*=0~~on~\partial D_1^*\cup\partial D_2^*,~~~~~v_3^*= 
\varphi^*~~on~\partial\Omega.~~~~~~~~~~\\ 
\end{aligned} 
\right. 
\end{equation} 
First we prove 
\begin{lem}\label{lm:Vi} 
There exist unique $v^*_i\in 
L^\infty(\widetilde \Omega^*)\cap 
C^0(\overline{ \widetilde \Omega^*}\setminus \{0\}) 
 \cap 
~C^2(\widetilde \Omega^*)$, $i=1,2,3$, 
which solve equations 
 $(\ref{eq:V1})$, $(\ref{eq:V2})$ and $(\ref{eq:V3})$ 
respectively. Moreover, 
$v_i^*\in 
C^1(\overline{ \widetilde \Omega^*}\setminus \{0\})$. 
\end{lem} 
\emph{Proof}:\hspace{0.2cm}The existence of solutions to the above 
equations can easily be obtained by Perron's method, 
see theorem 2.12 and lemma 2.13 in \cite{GT}. For reader$'$s convenience, we give 
below a simple proof of the 
 uniqueness. We only need to 
 prove that $0$ is the only solution 
in 
$ 
L^\infty(\widetilde \Omega^*)\cap 
C^0(\overline{ \widetilde \Omega^*}\setminus \{0\}) 
\cap 
~C^2(\widetilde \Omega^*)$ 
 to the 
following equation: 
\begin{equation} 
\left\{\begin{aligned} 
             \Delta w&=0 ~~~~~in~\widetilde\Omega^*, \\ 
             w&=0~~~~~on~\partial{\widetilde\Omega^*}\backslash\{0\}. 
\end{aligned} 
\right. 
\end{equation} 
Indeed, $\forall ~\varepsilon>0$, we have 
$$ |w(x)|\leq 
\frac{\varepsilon^{n-2}\|w\|_{L^\infty(\widetilde\Omega^*)}}{|x|^{n-2}},~~~on~\partial{(\widetilde\Omega^*\backslash 
B_\varepsilon)(0)}. 
$$ 
By the maximum principle,  $$ |w(x)|\leq 
\frac{\varepsilon^{n-2}\|w\|_{L^\infty(\widetilde\Omega^*)}}{|x|^{n-2}},~~~\forall~x\in\widetilde\Omega^*\backslash 
B_\varepsilon(0). 
$$ 
Thus $w\equiv 0$ in $\widetilde\Omega^*$. The additional regularity 
$v_i^*\in C^1(\overline{ \widetilde \Omega^*}\setminus \{0\})$ 
follows from standard elliptic estimates and the regularity of the 
$\partial D_i$ and $\partial \Omega$. $\hfill\square$ 
 
\begin{lem}\label{lem2.1}~For $i=1,2,3$, 
\begin{equation} 
v_{i\varepsilon}\longrightarrow v_i^* 
~~~\text{in}~~ C^2_{loc}(\widetilde \Omega^*), 
~~~~ 
\text{as}~~~\varepsilon\rightarrow 0, 
\label{a1} 
\end{equation} 
\begin{equation} 
\int_{\partial \Omega}\frac{\partial v_{i\varepsilon}}{\partial 
\nu}\longrightarrow\int_{\partial \Omega}\frac{\partial 
v^*_i}{\partial \nu}, ~~~~~\text{as}~~~ \varepsilon\rightarrow 0, 
~~~~i=1,2, 
\label{a2} 
\end{equation} 
\begin{equation} 
\int_{\partial D_{i\varepsilon}}\frac{\partial 
v_{3\varepsilon}}{\partial \nu}\longrightarrow\int_{\partial 
D^*_i}\frac{\partial v^*_3}{\partial \nu}, ~~~~~\text{as}~~~ 
\varepsilon\rightarrow 0. 
\label{a3} 
\end{equation} 
\end{lem} 
\emph{Proof}:~By the maximum principle, 
$\{\|v_{i\varepsilon}\|_{L^\infty}\}$ is bounded by a constant 
independent of $\varepsilon$. By the uniqueness part of Lemma 
\ref{lm:Vi}, we obtain (\ref{a1}) using standard elliptic estimates. 
By Lemma \ref{lm grad v3}, $\{\|\nabla v_{3\varepsilon}\|_{ 
L^\infty}\}$ is bounded by some constant independent of 
$\varepsilon$, so $\|\nabla v_3^*\|_{ L^\infty}<\infty$.  Estimate 
(\ref{a2}) and (\ref{a3}) follow from standard elliptic estimates. 
The proof is completed.$\hfill\square$\\\\ 
Similar to $Q_\varepsilon[\varphi_\varepsilon]$, we define 
\begin{equation} \label{eq:Q} 
Q^*[\varphi^*]:=\int_{\partial D^*_1}\frac{\partial v^*_3}{\partial 
\nu}\int_{\partial \Omega}\frac{\partial v^*_2}{\partial 
\nu}-\int_{\partial D^*_2}\frac{\partial v^*_3}{\partial 
\nu}\int_{\partial \Omega}\frac{\partial v^*_1}{\partial \nu}, 
\end{equation} 
then $Q^*: C^2(\partial\Omega)\mapsto \mathbb{R}$ is a linear 
functional. Let $Q_\varepsilon[\varphi_\varepsilon]$ and 
$Q^*[\varphi^*]$ be defined by equation (\ref{eq:Q_e}), 
(\ref{eq:Q}), then, by the above lemmas, 
$$Q_\varepsilon[\varphi_\varepsilon]\longrightarrow Q^*[\varphi^*], 
~~~~~\text{as}~~~ \varepsilon\rightarrow 0.$$ 
 
\begin{cor}\label{cor:Q} 
If  $\varphi^*\in C^2(\partial\Omega)$ satisfies 
 $Q^*[\varphi^*]\neq0$, then $|Q_\varepsilon[\varphi_\varepsilon]|\geq \frac{1}{C}$, for 
some positive constant $C$ which is independent of $\varepsilon$. 
\end{cor} 
 
In the following we give some examples to show that, in general, the 
rates of the lower bounds established in Theorem \ref{thm:lowbdd} 
are optimal.\\ 
 
Let $\Omega\subset\mathbb{R}^n$, $n\geq 2$, be a bounded open set 
with $C^{2,\alpha}$ boundary, $0<\alpha<1$, which is symmetric with 
respect to $x_1$-variable, i.e., $(x_1,x')\in \Omega$ if and only if 
$(-x_1,x')\in \Omega$, where $x'=(x_2,\cdots,x_n)$. 
 
Let $D_1^*$ be a strictly convex bounded open set in $\{(x_1,x')\in 
\mathbb{R}^n |x_1<0\}$ with $C^{2,\alpha}$ boundary, $0<\alpha<1$, 
satisfying $0\in\partial D_1^*$ and $\overline 
{D_1^*}\subset\Omega$. Set $D_2^*=\{(x_1,x')\in \mathbb{R}^n 
|(-x_1,x')\in D_1^*\}$. 
 
Let $\varphi\in C^2(\partial\Omega)\backslash\{0\}$ satisfy 
\begin{equation} \label{eq:varphi odd} 
\varphi_{odd}(x_1,x'):=\frac{1}{2}\big[\varphi(x_1,x')-\varphi(-x_1,x')\big]\leq 
0 ~(\text{or} \geq0), 
\end{equation} 
on $(\partial \Omega)^+:=\{(x_1,x')\in \partial\Omega|x_1>0\}$. 
 
For $\varepsilon>0$ sufficiently small, let 
\begin{equation} 
\begin{split} 
D_{1\varepsilon}&:=\big\{(x_1,x')\in 
\Omega\big|(x_1+\frac{\varepsilon}{2},x')\in D_1^*\big\},\\ 
D_{2\varepsilon}&:=\big\{(x_1,x')\in 
\Omega\big|(x_1-\frac{\varepsilon}{2},x')\in D_2^*\big\},\\ 
\varphi_\varepsilon&:=\varphi. 
\end{split} 
\nonumber 
\end{equation} 
 
\begin{prop} \label{prop:varphi odd} 
Under the above assumptions, we have $|Q_\varepsilon[\varphi]|\geq 
\frac{1}{C}$, for some positive constant $C$ independent of 
$\varepsilon$. Consequently, 
\begin{equation} \label{eq:lower epsilon} 
\begin{aligned} 
           &\|\nabla u_\varepsilon\|_{L^{\infty}(\widetilde{\Omega})}\geq 
           \frac{1}{C\sqrt{\varepsilon}}~~~~~~~~~for \hspace{0.2cm} n=2,\\ 
           &\|\nabla u_\varepsilon\|_{L^{\infty}(\widetilde{\Omega})}\geq 
           \frac{1}{C\varepsilon|\ln{\varepsilon}|}~~~~~for \hspace{0.2cm} n=3,\\ 
           &\|\nabla u_\varepsilon\|_{L^{\infty}(\widetilde{\Omega})}\geq 
           \frac{1}{C\varepsilon}~~~~~~~~~~~for \hspace{0.2cm} n \geq 
           4, 
\end{aligned} 
\end{equation} 
where $u_\varepsilon$ is the solution to equation (\ref{eq:k 
+infty}). 
\end{prop} 
 
The above proposition can be easily obtained by the following lemma 
which gives a necessary and sufficient condition instead of 
condition (\ref{eq:varphi odd}) on $\varphi$ for the lower bounds 
(\ref{eq:lower epsilon}) to hold. 
 
Let 
\begin{equation}\label{eq:v3 odd} 
(v^*_3)_{odd}(x_1,x'):=\frac{1}{2}\big[v^*_3(x_1,x')-v^*_3(-x_1,x')\big], 
\end{equation} 
we have 
\begin{lem}\label{lm:W} 
Under the same hypotheses in Proposition \ref{prop:varphi odd} 
except for the condition (\ref{eq:varphi odd}), let 
$Q_\varepsilon[\varphi]$ and $(v^*_3)_{odd}(x)$ be defined by 
equation (\ref{eq:Q_e}) and (\ref{eq:v3 odd}), then the following 
statements are equivalent: 
\begin{enumerate} 
\item For some positive constant $C$ independent of 
$\varepsilon$, we have $|Q_\varepsilon[\varphi]|\geq \frac{1}{C}$, 
\item $\int_{\partial D^*_2}\frac{\partial (v^*_3)_{odd}}{\partial 
\nu}\neq0$. 
\end{enumerate} 
\end{lem} 
\emph{Proof}:\hspace{0.2cm} By symmetry, the strong maximum 
principle and the Hopf Lemma, we can easily obtain 
$$\int_{\partial 
\Omega}\frac{\partial v^*_1}{\partial \nu}=\int_{\partial 
\Omega}\frac{\partial v^*_2}{\partial \nu}<0.$$ Then 
\begin{equation}\begin{split} 
Q^*[\varphi]&=\int_{\partial \Omega}\frac{\partial v^*_1}{\partial 
\nu}\Big(\int_{\partial D^*_1}\frac{\partial v^*_3}{\partial 
\nu}-\int_{\partial D^*_2}\frac{\partial v^*_3}{\partial 
\nu}\Big)\\ 
&=\int_{\partial \Omega}\frac{\partial v^*_1}{\partial 
\nu}\Big(\int_{\partial D^*_1}\frac{\partial (v^*_3)_{odd}}{\partial 
\nu}-\int_{\partial D^*_2}\frac{\partial (v^*_3)_{odd}}{\partial 
\nu}\Big)\\ 
&=-2\int_{\partial \Omega}\frac{\partial v^*_1}{\partial 
\nu}\int_{\partial D^*_2}\frac{\partial (v^*_3)_{odd}}{\partial 
\nu}. 
\end{split} 
\nonumber 
\end{equation} 
Hence, $Q^*[\varphi]\neq 0~$ if and only if $~\int_{\partial 
D^*_2}\frac{\partial 
(v^*_3)_{odd}}{\partial \nu}\neq0$. Then by Corollary \ref{cor:Q}, we complete the 
proof.$\hfill\square$\\\\ 
\emph{Proof of Proposition \ref{prop:varphi odd}}:\hspace{0.2cm} 
Note that $(v^*_3)_{odd}(0,x')=0$ by symmetry, and $(v^*_3)_{odd}$ 
is harmonic with $(v^*_3)_{odd}=\varphi_{odd}\leq 0~(\text{or} 
\geq0)$ but not identically 0 on $(\partial\Omega)^+$. Now by using 
the strong maximum principle and the Hopf Lemma, it is clear that 
$\int_{\partial D^*_2}\frac{\partial (v^*_3)_{odd}}{\partial 
\nu}\neq0$, Hence, by Lemma \ref{lm:W} and Theorem \ref{thm:lowbdd}, 
we are done.$\hfill\square$\\ 
 
\begin{rem} 
If $\varphi=\sum_{i=1}^{n}b_ix_i$ with $b_i\in \mathbb{R}$ and 
$b_1\neq 0$, then by Proposition \ref{prop:varphi odd} we have 
$|Q_\varepsilon[\varphi]|\geq \frac{1}{C}$. Therefore, by Theorem 
\ref{thm:upbdd} and \ref{thm:lowbdd}, the blow-up rates of $\|\nabla 
u\|_{L^\infty(\widetilde\Omega)}$ are $\varepsilon^{-1/2}$ in in 
dimension $n=2$, $(\varepsilon|\ln\varepsilon|)^{-1}$ in dimension 
$n=3$ and $\varepsilon^{-1}$ in dimension $n\geq4$. 
\end{rem} 
 
Now instead of in a bounded set $\Omega$, we consider in 
$\mathbb{R}^n$: 
\begin{equation} \label{kang eqn} 
\left\{\begin{aligned} 
             &\Delta u_\varepsilon=0~~~~~~~~~~~~~~~~in~\mathbb{R}^n \backslash 
             \overline{D_{1\varepsilon}\cup D_{2\varepsilon}}, \\ 
                 &u_\varepsilon|_{+}=u_\varepsilon|_{-}~~~~~~~~~~~~on~\partial D_{1\varepsilon}\cup \partial D_{2\varepsilon}, \\ 
           &\nabla u_\varepsilon\equiv 0~~~~~~~~~~~~~~~~in~ D_{1\varepsilon}\cup D_{2\varepsilon}, \\ 
           &\int_{\partial D_{i\varepsilon}}\frac{\partial u_\varepsilon}{\partial\nu}\Big|_{+}=0 
           ~~~~~~~~~~(i=1,2), \\ 
                 &\limsup_{|x|\to \infty}|x|^{n-1}|u_\varepsilon(x)-H(x)|< \infty, 
\end{aligned} 
\right. 
\end{equation} 
where $H(x)$ is a given entire harmonic function in $\mathbb{R}^n$. 
 
we have the following result regarding the lower bound for $|\nabla 
u_\varepsilon|$: 
\begin{prop}\label{p.H} 
With the same assumptions on $D_{1\varepsilon}$ and 
$D_{2\varepsilon}$ as in Proposition \ref{prop:varphi odd}, and let 
$H(x)$ be an entire harmonic function in $\mathbb{R}^n$ satisfying 
$H_{odd}(x_1,x'):=\frac{1}{2}\big[H(x_1,x')-H(-x_1,x')\big]< 0 
~(\text{or} > 0)$ on $\mathbb{R}^n_+:=\{(x_1,x')\in 
\mathbb{R}^n|x_1>0\}$, then for some positive constant $C$ 
independent of $\varepsilon$, we have 
\begin{equation} \label{eq:Kang lower epsilon} 
\begin{aligned} 
           &\|\nabla u_\varepsilon\|_{L^{\infty}(\mathbb{R}^n \backslash 
             \overline{D_{1\varepsilon}\cup D_{2\varepsilon}})}\geq 
           \frac{1}{C\sqrt{\varepsilon}}~~~~~~~~~for \hspace{0.2cm} n=2,\\ 
           &\|\nabla u_\varepsilon\|_{L^{\infty}(\mathbb{R}^n \backslash 
             \overline{D_{1\varepsilon}\cup D_{2\varepsilon}})}\geq 
           \frac{1}{C\varepsilon|\ln{\varepsilon}|}~~~~~for \hspace{0.2cm} n=3,\\ 
           &\|\nabla u_\varepsilon\|_{L^{\infty}(\mathbb{R}^n \backslash 
             \overline{D_{1\varepsilon}\cup D_{2\varepsilon}})}\geq 
           \frac{1}{C\varepsilon}~~~~~~~~~~~for \hspace{0.2cm} n \geq 
           4, 
\end{aligned} 
\end{equation} 
where $u_\varepsilon$ is the solution to equation (\ref{kang eqn}). 
\end{prop} 
\emph{Proof}: Step 1. First, we show that there exists a positive 
constant $C$ independent of $\varepsilon$, such that for any small $ 
\varepsilon>0$, 
\begin{equation}\label{eq:uniform bd w} 
|x|^{n-1}|u_\varepsilon(x)-H(x)|\leq C,\qquad \forall\ x\in 
\mathbb{R}^n \backslash\overline{D_{1\varepsilon}\cup 
D_{2\varepsilon}}. 
\end{equation} 
(i) For any bounded open set  $U\subset\mathbb{R}^n$ with $C^1$ 
boundary $\partial U$ satisfying $\partial U\cap 
\overline{D_{1\varepsilon}\cup D_{2\varepsilon}}=\emptyset$, we 
have, in view of the first and the fourth lines in (\ref{kang eqn}), 
\begin{equation} \label{eq:integral normal u=0} 
\int_{\partial U}\frac{\partial u_\varepsilon}{\partial\nu} = \int_{ 
U\setminus \overline{D_{1\varepsilon}\cup D_{2\varepsilon}} }\Delta 
u_\varepsilon =0. 
\end{equation} 
(ii) We show that there exists a positive constant $M$ independent 
of $\varepsilon$, such that 
$$\|u_\varepsilon 
-H\|_{L^\infty(\mathbb{R}^n \backslash 
             \overline{D_{1\varepsilon}\cup D_{2\varepsilon}})}\leq M, ~~~~~\forall ~\text{small}~\varepsilon>0.$$ 
We only need to prove \begin{equation} \|u_\varepsilon 
-H\|_{L^\infty(\mathbb{R}^n \backslash 
             \overline{D_{1\varepsilon}\cup D_{2\varepsilon}})}\leq 
             \sum_{i=1}^2 (\max_{\overline D_{i\varepsilon} }H- 
             \min_{\overline D_{i\varepsilon} }H). 
             \label{B1} 
             \end{equation} 
Since $\nabla u_\varepsilon=0$ in $D_{1\varepsilon}\cup 
D_{2\varepsilon}$, $u_\varepsilon$ is constant on each 
$D_{i\varepsilon}$, denoted as $C_i(\varepsilon)$.  We know that 
\begin{equation}\lim_{|x|\to \infty}(u_\varepsilon(x)-H(x))=0, 
\label{B2} \end{equation} 
 and \begin{equation} C_i(\varepsilon)-\max_{ \overline 
D_{i\varepsilon} }H \le u_\varepsilon-H\le  C_i(\varepsilon)-\min_{ 
\overline D_{i\varepsilon} }H, \quad \text{on}\ D_{i\varepsilon},~~~ 
i=1,2. \label{B3} \end{equation} 
 If (\ref{B1}) did not hold, say, $$ 
\sup_{\mathbb{R}^n}(u_\varepsilon-H) 
> 
 \sum_{i=1}^2 (\max_{\overline D_{i\varepsilon} }H- 
             \min_{\overline D_{i\varepsilon} }H), 
             $$ 
             then, because of (\ref{B2}) and (\ref{B3}),  there would exist $\displaystyle{ 
             0<a<\sup_{\mathbb{R}^n}(u_\varepsilon-H) 
             }$ such that 
             $U:= 
 \{x\in \mathbb{R}^n\ |\ (u_\varepsilon-H)(x)>a\} 
 \neq \emptyset $ satisfies 
 $\partial U\cap 
\overline{D_{1\varepsilon}\cup D_{2\varepsilon}}=\emptyset$.  We may 
assume, by the Sard theorem, that $a$ is a regular value of 
$u_\varepsilon-H$, and therefore $\partial U$ is $C^1$.  By the Hopf 
lemma, $\displaystyle{ \frac{  \partial (u_\varepsilon-H) }{\partial 
\nu}<0}$ on $\partial U$, and therefore 
$$ 
\int_{\partial U}\frac{  \partial (u_\varepsilon-H) }{\partial 
\nu}<0. 
$$ 
On the other hand, using (\ref{eq:integral normal u=0}) and the 
harmonicity of $H$ in $U$, we have 
$$ 
\int_{\partial U}\frac{  \partial (u_\varepsilon-H) }{\partial \nu} 
=-\int_{\partial U}\frac{  \partial H }{\partial \nu} =-\int_U 
\Delta H=0. 
$$ 
A contradiction.\\\\ 
(iii) Consider $w_\varepsilon(x):=u_\varepsilon(x)-H(x)$. Fix a 
constant  $R_0>0$, independent of $\varepsilon$, such that 
$D_1^*\cup D_2^*\subset B_{ R_0/2}(0)$, and  let 
$$\widetilde{w_\varepsilon}(y):=\frac{1}{|y|^{n-2}}w_\varepsilon\Big(\frac{y}{|y|^2}\Big), ~~~0<|y|<\frac{1}{R_0}.$$ 
Then $\widetilde{w_\varepsilon}$ is harmonic in 
$B_{1/R_0}\backslash\{0\}$. By the last line of (\ref{kang eqn}), 
there exists a positive constant $C(\varepsilon)$  such that 
$$|\widetilde{w_\varepsilon}(y)|\leq C(\varepsilon)|y|, ~~~0<|y|<\frac{1}{R_0}.$$ 
Therefore, $\Delta \widetilde{w_\varepsilon}=0$ in $B_{1/R_0}$ and 
$\widetilde{w_\varepsilon}(0)=0$. By (ii), we have 
$|\widetilde{w_\varepsilon}|\leq C$, on $\partial B_{1/R_0}$, for 
some positive constant $C$ independent of $\varepsilon$. Hence, 
$|\widetilde{w_\varepsilon}|\leq C, 
~|\nabla\widetilde{w_\varepsilon}|\leq C$ in $B_{1/(2R_0)}$, then 
$$|\widetilde{w_\varepsilon}(y)|\leq C|y|, ~~~|y|<\frac{1}{2R_0}.$$ 
Therefore, also using (ii), (\ref{eq:uniform bd w}) holds.\\\\ 
Step 2. For $R>R_0$, let  $\Omega=B_R(0)$. Let 
$\varphi_\varepsilon:=u_\varepsilon|_{\partial \Omega}$, then by 
Corollary \ref{cor:Q} and Theorem \ref{thm:lowbdd} it is enough to 
show, for some $R$,  that $Q^*[\varphi^*]\neq 0$, where $\varphi^*$ 
is defined at the beginning of this section. By symmetry, we have 
$$Q^*[\varphi^*]=\int_{\partial \Omega}\frac{\partial v^*_1}{\partial \nu}\Big(\int_{\partial D^*_1}\frac{\partial v^*_3}{\partial 
\nu}-\int_{\partial D^*_2}\frac{\partial v^*_3}{\partial 
\nu}\Big).$$ 
 
Without loss of generality, we may assume $H_{odd}(x)> 0$ on 
$\mathbb{R}^n_+$.  Recall that $v_3^*$ is the solution of 
(\ref{eq:V3}) with boundary data $\varphi^*$.  In the following we 
use notation $(v_3^*)_h$ to denote the the solution of (\ref{eq:V3}) 
with boundary data $h$.  Since $Q^*[\varphi^*]$ is linear on 
$\varphi^*$ and by symmetry $Q^*[H_{even}]=H[ \varphi^*_{even}]=0$, 
where 
$H_{even}(x):=H(x)-H_{odd}(x)=\frac{1}{2}\big[H(x_1,x')+H(-x_1,x')\big]$ 
and similar for $\varphi^*_{even}$, we may assume $H(x)=H_{odd}(x)$.

Now consider $w(x)=H(x)-(v_3^*)_H(x)$. Then $w(x)$ is harmonic in 
$\widetilde\Omega^*$ which is defined at the beginning of this 
section. By symmetry, $w(-x_1,x')=-w(x_1,x')$, $w(x)=H(x)$ on 
$\partial D_1^*\cup\partial D_2^*$ and $w(x)=0$ on $\partial 
\Omega$. Therefore, 
$$-2\int_{\partial {D_2^*}} 
H\frac{\partial w}{\partial \nu}=\int_{\widetilde\Omega^*} 
w(x)\Delta w(x)+ \int_{\widetilde\Omega^*} |\nabla w|^2= 
\int_{\widetilde\Omega^*} |\nabla w|^2\ge 0.$$ On the other hand, 
$(v_3^*)_H=0$ on $\partial D_2^*$, $(v_3^*)_H>0$ on $(\partial 
\Omega)^+$ and, by the oddness of $(v_3^*)_H$, $(v_3^*)_H=0$ on 
$\{(x_1, x')\ |\ x_1=0\}$.  Thus, by the maximum principle and the 
strong maximum principle, $(v_3^*)_H>0$ in $\widetilde\Omega^*$ and 
in turn, using the Hopf lemma, $\displaystyle{ 
 \frac{\partial 
(v_3^*)_H}{\partial \nu}>0 }$ on $\partial {D_2^*}$. 
 Hence, using the harmonicity of $H$, 
\begin{equation} 
\begin{split} 
\max_{\partial D_2^*}H\int_{\partial {D_2^*}} \frac{\partial 
(v_3^*)_H}{\partial \nu} &\geq\int_{\partial {D_2^*}} H\frac{\partial 
(v_3^*)_H}{\partial \nu}\geq \int_{\partial {D_2^*}} H\frac{\partial H}{\partial 
\nu}-\int_{\partial {D_2^*}} H\frac{\partial w}{\partial 
\nu}\\ 
&\geq\int_{D_2^*} |\nabla H|^2\geq \frac{1}{C}, 
\end{split} 
\nonumber 
\end{equation} 
Therefore, $$\int_{\partial {D_2^*}} \frac{\partial 
(v_3^*)_H}{\partial \nu}\geq\frac{1}{C},$$ for positive constant $C$ 
independent of $R$. 
 
For $s_\varepsilon:=\varphi_\varepsilon-H$ on $\partial \Omega$, by step 1, there exists a constant $C>0$ which is independent of $\varepsilon$ and $R$, 
 suth that $\|s_\varepsilon\|_{L^\infty(\partial\Omega)}\leq CR^{1-n}$. By Remark \ref{rem1.1}, we have $\|\nabla 
(v_3^*)_{s^*}\|_{L^\infty(\partial D_1^*\cup \partial D_2^*)}\leq 
C\|s^*\|_{L^\infty(\partial\Omega)}$, thus, 
$$\Big|\int_{\partial {D_i^*}} \frac{\partial (v_3^*)_{s^*}}{\partial \nu}\Big|\leq 
C\int_{\partial {D_i^*}} \|s^*\|_{L^\infty(\partial\Omega)}\leq CR^{1-n},$$ for some positive constant $C$ independent of $\varepsilon$ and $R$. 
 
Therefore, for large enough $R$, 
$$\int_{\partial {D_2^*}} \frac{\partial 
(v_3^*)_{\varphi^*}}{\partial \nu}=\int_{\partial {D_2^*}} 
\frac{\partial (v_3^*)_H}{\partial \nu}+\int_{\partial {D_2^*}} 
\frac{\partial (v_3^*)_{s^*}}{\partial \nu}\geq\frac{1}{C}\neq 0.$$ 
 
It is also clear that $\int_{\partial \Omega}\frac{\partial 
v^*_1}{\partial \nu}<0$, Thus, 
$$Q^*[\varphi^*]=-2\int_{\partial \Omega}\frac{\partial v^*_1}{\partial \nu} 
\int_{\partial D^*_2}\frac{\partial (v^*_3)_{\varphi^*}}{\partial 
\nu}\neq 0.$$ This proof is completed.$\hfill\square$ 
 
\begin{rem} 
In $\mathbb{R}^2$, when $D_{1\varepsilon}$ and $D_{2\varepsilon}$ 
are identical balls of radius 1, the estimate (\ref{eq:Kang lower 
epsilon}) was established in \cite{AKLLL} under a weaker assumption 
$\partial_{x_1} H(0)\neq 0.$\\\\ 
\end{rem} 
 
%%%%%%%%%%%%%%%%%%%%%%%%%%%%%%%%%%%%%%%%%%%%%%%%%%%%%%%%%%%%%%%%%%%%%%%% 
% 
%      3. Proof of Theorem for general coefficients 
% 
%%%%%%%%%%%%%%%%%%%%%%%%%%%%%%%%%%%%%%%%%%%%%%%%%%%%%%%%%%%%%%%%%%%%%%%% 
 
\section {Proof of Theorem \ref{thm:upbdd general} and \ref{thm:lowbdd general}} 
 
In the introduction, similar to the harmonic case, we still 
decompose $u=C_1V_1+C_2V_2+V_3$ as in (\ref{eq:decomp general}).\\\\ 
Proposition \ref{prop:bd C_1-C_2} holds since Lemma \ref{lm grad 
v_1,2}$-$\ref{lm grad v3} hold for $V_1$, $V_2$, $V_3$ defined by 
(\ref{eq:V1 general})$-$(\ref{eq:V3 general}) and $\rho\in 
C^2(\widetilde{\Omega})$ which is the solution to: 
\begin{equation} 
        \left\{ \begin{aligned} 
           &\partial_{x_j}\Big(a_2^{ij}(x)~\partial_{x_i}{\rho}\Big)=0~~~~~~~~in~\widetilde\Omega,\\ 
           &\rho=0 ~~on~\partial D_1\cup\partial D_2,~~~~\rho=1 ~~on~\partial\Omega.\\ 
          \end{aligned} 
\right. \nonumber 
\end{equation} 
The proofs are essentially the same.\\ 
 
Now we start to estimate $|C_1-C_2|$. By the decomposition formula 
(\ref{eq:decomp general}), instead of (\ref{eq:ab}), we denote 
\begin{equation} \label{eq:ab general} 
\begin{split} 
a_{lm}&=\int_{\partial 
D_l}a_2^{ij}(x)~\partial_{x_i}{V_m}~\nu_{j}~~(l,m=1,2),\\ 
b_l&=\int_{\partial 
D_l}a_2^{ij}(x)~\partial_{x_i}{V_3}~\nu_{j}~~~(l=1,2). 
\end{split} 
\end{equation} 
Then Lemma \ref{lm:ab} and (\ref{eq:c1c2})$-$(\ref{eq:c1-c2}) still 
hold for $a_{lm}$ and $b_l$ defined 
above.\\\\ 
In fact, to prove Lemma \ref{lm:ab} with general coefficients, we 
only need to change $\frac{\partial *}{\partial \nu}$ to 
$a_2^{ij}(x)~\partial_{x_i}{*}~\nu_{j}$, change $\Delta *$ in 
$\partial_{x_j}\Big(a_2^{ij}(x)~\partial_{x_i}{*}\Big)$ and change 
$v_1$, $v_2$, $v_3$ in $V_1$, $V_2$, $V_3$, respectively, in the 
original proof of Lemma \ref{lm:ab}. For instance, 
(\ref{eq:a12=a21}) is changed to 
\begin{equation} \label{eq:a12=a21 general} 
\begin{split} 
0&=\int_{\widetilde\Omega}\partial_{x_j}\Big(a_2^{ij}(x)~\partial_{x_i}{V_1}\Big)\cdot 
V_2 - 
\int_{\widetilde\Omega}\partial_{x_j}\Big(a_2^{ij}(x)~\partial_{x_i}{V_2}\Big)\cdot 
V_1\\ 
&= -\int_{\partial D_2}a_2^{ij}(x)~\partial_{x_i}{V_1}~\nu_{j}\cdot 
1 + \int_{\partial D_1}a_2^{ij}(x)~\partial_{x_i}{V_2}~\nu_{j}\cdot 
1\\ 
&=-a_{21}+a_{12}. 
\end{split} 
\end{equation} 
Therefore, to estimate $|C_1-C_2|$, it is equivalent to estimating 
$|a_{11}-\alpha a_{12}|$ and $|b_1-\alpha b_2|$.\\ 
 
For $|a_{11}-\alpha a_{12}|$, Lemma \ref{lm: bd aii n 2}$-$\ref{lm: 
bd aii n 4+} still hold for $a_{ll}(l=1,2)$ defined by (\ref{eq:ab 
general}). The proof is quite similar and the only thing which needs 
to be shown is the following: 
\begin{equation} 
\begin{split} 
0&=\int_{\widetilde\Omega}\partial_{x_j}\Big(a_2^{ij}(x)~\partial_{x_i}{V_1}\Big)\cdot 
V_1\\ 
&=-\int_{\widetilde\Omega}a_2^{ij}(x)~\partial_{x_i}{V_1}\partial_{x_j}{V_1} 
-\int_{\partial D_1}a_2^{ij}(x)~\partial_{x_i}{V_1}~\nu_{j}\cdot 
1\\ 
&=-\int_{\widetilde\Omega}a_2^{ij}(x)~\partial_{x_i}{V_1}\partial_{x_j}{V_1}-a_{11}, 
\end{split} 
\nonumber 
\end{equation} 
i.e. 
$$a_{11}=-\int_{\widetilde\Omega}a_2^{ij}(x)~\partial_{x_i}{V_1}\partial_{x_j}{V_1}.$$ 
Then by the uniform ellipticity of $a_2^{ij}(x)$ and the harmonicity 
of $v_1$, 
$$ 
|a_{11}|\geq\lambda\int_{\widetilde\Omega}|\nabla 
V_1|^2\geq\lambda\int_{\widetilde\Omega}|\nabla v_1|^2, 
$$ 
and 
$$|a_{11}|\leq\int_{\widetilde\Omega}a_2^{ij}(x)~\partial_{x_i}{w}\partial_{x_j}{w} 
\leq\Lambda\int_{\widetilde\Omega}|\nabla 
w|^2\leq\Lambda\int_{\widetilde\Omega\cap O_{r/2}}|\nabla \overline 
w|^2+C,$$ where $w$ is defined in the proof of Lemma \ref{lm: bd aii 
n 2} with the same boundary data of $V_1$ and $\overline w$ is 
defined by (\ref{eq:w n=2}) and (\ref{eq:w n>2}).\\\\ 
Thus, Lemma \ref{lm: bd aii n 2}$-$\ref{lm: bd aii n 4+} follow by 
the same computations. Then Lemma \ref{lm:bd alpha} and Proposition 
\ref{prop: bd a11-alpha a22} hold with 
the same proofs.\\ 
 
For $|b_{1}-\alpha b_{2}|$, Proposition \ref{prop:bd b1-alpha b2} 
also holds for $b_l(l=1,2)$ defined by (\ref{eq:ab general}) and 
$Q_\varepsilon[\varphi]$ defined by (\ref{eq:Q_e general}). The 
proof is the same after changing $\frac{\partial *}{\partial \nu}$ 
to 
$a_2^{ij}(x)~\partial_{x_i}{*}~\nu_{j}$.\\ 
 
Combining the above propositions, we obtain our theorems.\\\\

%%%%%%%%%%%%%%%%%%%%%%%%%%%%%%%%%%%%%%%%%%%%%%%%%%%%%%%%%%%%%%%%%%%%%%%% 
% 
%      4. Appendix A 
% 
%%%%%%%%%%%%%%%%%%%%%%%%%%%%%%%%%%%%%%%%%%%%%%%%%%%%%%%%%%%%%%%%%%%%%%%% 
 
\section{Appendix} 
{\large\bf Some elementary results for the conductivity problem} 
 
\vspace{0.3cm} Assume that in $\mathbb{R}^n$, $\Omega$ and $\omega$ 
are  bounded open sets with  $C^{2,\alpha}$ boundaries, 
$0<\alpha<1$, satisfying 
$$\overline\omega=\bigcup_{s=1}^{m}\overline\omega_s\subset\Omega,$$ 
where $\{\omega_s\}$ are 
 connected components of $\omega$. 
Clearly, $m<\infty$ and $\omega_s$ is open for all $1\le s\le \omega$. 
 Given $\varphi\in 
C^2(\partial\Omega)$, the conductivity problem we consider is the 
following transmission problem with Dirichlet boundary condition: 
\begin{equation} \label{diveq:k finite A} 
\left\{ \begin{aligned} 
           \partial_{x_j}\Big\{\Big[\big(ka_1^{ij}(x)-a_2^{ij}(x)\big) 
\chi_\omega+a_2^{ij}(x)\Big]\partial_{x_i}u_k\Big\}&=0~~~in~\Omega, \\ 
           u_k=\varphi ~~~~~~~~~~~~~~~~~~~&~~~~~~~~on~\partial\Omega, 
          \end{aligned} 
\right. 
\end{equation} 
where 
%\begin{equation} 
%a_k(x)=\left\{ \begin{aligned} 
%           &kA_1(x)=k\big(a_1^{ij}(x)\big)~~~~in~\omega, \\ 
%           &A_2(x)=\big(a_2^{ij}(x)\big)~~~~~~~in~\Omega\backslash\overline{\omega}, 
%          \end{aligned} 
%\right. \nonumber 
%\end{equation} 
$k=1,2,3,\cdots$, and $\chi_\omega$ is the characteristic function of 
$\omega$. 
 
The $n\times n$ matrixes 
$A_1(x):=\big(a_1^{ij}(x)\big)~\text{in}~\omega,~ 
A_2(x):=\big(a_2^{ij}(x)\big)~\text{in}~\Omega\backslash\overline{\omega}$ 
are symmetric and $\exists$ a constant $\Lambda\geq\lambda>0$ such 
that 
$$~\lambda|\xi|^2\leq a_1^{ij}(x)\xi_i\xi_j\leq \Lambda|\xi|^2~(\forall x\in\omega), 
~~~~\lambda|\xi|^2\leq a_2^{ij}(x)\xi_i\xi_j\leq 
\Lambda|\xi|^2~(\forall x\in\Omega\backslash\omega)$$ for all 
$\xi\in\mathbb{R}^n$ and $a_1^{ij}(x)\in C^2(\overline \omega), 
~a_2^{ij}(x) \in C^2(\overline \Omega\backslash\omega)$. 
 
\vspace{0.2cm} Equation (\ref{diveq:k finite A}) can be rewritten in 
the following form to emphasize the transmission condition on 
$\partial\omega$: 
\begin{equation} \label{eq:k finite A} 
\left\{ \begin{aligned} 
           &\partial_{x_j}\Big(a_1^{ij}(x)~\partial_{x_i}{u_k}\Big)=0~~~~~~~~~~~~~~~~~~~~in~\omega,\\ 
           &\partial_{x_j}\Big(a_2^{ij}(x)~\partial_{x_i}{u_k}\Big)=0~~~~~~~~~~~~~~~~~~~~in~\Omega 
           \backslash\overline\omega,\\ 
           &u_k|_{+}=u_k|_{-},~~~~~~~~~~~~~~~~~~~~~~~~~~~~~~~on~\partial\omega, \\ 
           &a_2^{ij}(x)\partial_{x_i}{u_k}\nu_j\big|_{+}=ka_1^{ij}(x)\partial_{x_i}{u_k}\nu_j\big|_{-}~~~on~\partial\omega, \\ 
           &u_k=\varphi~~~~~~~~~~~~~~~~~~~~~~~~~~~~~~~~~~~~~~~on~\partial\Omega. 
          \end{aligned} 
\right. 
%\nonumber 
\end{equation} 
Here and throughout this paper $\nu$ is the outward unit normal and 
the subscript $\pm$ indicates the limit from outside and inside the 
domain, respectively.\\\\ 
We list the following results which are well known and omit the 
proofs. 
\begin{thm} \label{thm:regular k finite} If $u_k\in H^1(\Omega)$ is a solution of equation (\ref{diveq:k finite 
A}), then $u_k\in C^1(\overline{\Omega\backslash\omega})\cap 
C^1(\overline\omega)$ and satisfies equation (\ref{eq:k finite A}). 
 
If $u_k\in C^1(\overline{\Omega\backslash\omega})\cap 
C^1(\overline\omega)$ is a solution of equation (\ref{eq:k finite 
A}), then $u_k\in H^1(\Omega)$ and satisfies equation (\ref{diveq:k 
finite A}). 
\end{thm} 
\begin{thm} \label{thm:unique k finite} There exists at most one solution $u_k\in H^1(\Omega)$ 
to equation (\ref{diveq:k finite A}). 
\end{thm} 
 
The existence of the solution can be obtained by using the 
variational method. For every $k$, we define the energy functional 
\begin{equation} \label{vm:k finite A} 
\begin{aligned} 
I_{k}[v]:&=\frac{k}{2}\int_{\omega}a_1^{ij}(x)\partial_{x_i}{v}\partial_{x_j}{v} 
+ 
\frac{1}{2}\int_{\Omega\backslash\overline{\omega}}a_2^{ij}(x)\partial_{x_i}{v}\partial_{x_j}{v}, 
\end{aligned} 
%\nonumber 
\end{equation} 
where $v$ belongs to the set 
$$H^1_{\varphi}(\Omega):=\{v\in H^1(\Omega)| ~v=\varphi ~~on~ \partial\Omega\}.$$ 
\begin{thm} \label{lem:exist k finite A} 
For every $k$, there exists a minimizer $u_k\in H^1(\Omega)$ 
satisfying 
$$I_{k}[u_k]=\min_{v\in H^1_{\varphi}(\Omega)}I_{k}[v].$$ 
Moreover, $u_k\in H^1(\Omega)$ is a solution of equation 
(\ref{diveq:k finite A}). 
\end{thm} 
 
Comparing equation (\ref{eq:k finite A}), when $k=+\infty$, the 
perfectly conducting problem turns out to be: 
\begin{equation} \label{eq:k +infty A} 
\left\{ \begin{aligned} 
           &\partial_{x_j}\Big(a_2^{ij}(x)~\partial_{x_i}{u}\Big)=0~~~~~~~~in~\Omega 
           \backslash\overline\omega,\\ 
           &u|_{+}=u|_{-}~~~~~~~~~~~~~~~~~~~~~~on~\partial\omega, \\ 
           &\nabla u=0~~~~~~~~~~~~~~~~~~~~~~~~~in~\omega, \\ 
           &\int_{\partial\omega_s}a_2^{ij}(x)\partial_{x_i}{u}\nu_j\big|_{+}=0 ~~~~(s=1,2,\cdots,m),\\ 
           &u=\varphi~~~~~~~~~~~~~~~~~~~~~~~~~~~on~\partial\Omega. 
          \end{aligned} 
\right. 
%\nonumber 
\end{equation} 
We also have similar results: 
\begin{thm} \label{thm:regular k +infty} If $u\in H^1(\Omega)$ satisfies 
equation (\ref{eq:k +infty A}) except for the fourth line, then 
$u\in C^1(\overline{\Omega\backslash\omega})\cap 
C^1(\overline\omega)$. 
\end{thm} 
\emph{Proof}:\hspace{0.2cm} By the third line of equation (\ref{eq:k 
+infty A}), we have $u\equiv \text{const}$ on each component of 
$\omega$, so $u\equiv \text{const}$ on each component of 
$\partial\omega$. Thus $u\equiv \text{const}$ on each component of 
$\partial(\Omega\backslash\overline\omega)$. 
 
Since $u\in H^1(\Omega)$ satisfies 
$\partial_{x_i}\Big(a_2^{ij}(x)~\partial_{x_i}{u_k}\Big)=0$ in 
$\Omega\backslash\overline\omega$, $u|_{\partial\Omega}=\varphi\in 
C^2(\partial\Omega)$ and $u\equiv \text{const}$ on each component of 
$\partial(\Omega\backslash\overline\omega)$, by the elliptic 
regularity theory, we have $u\in 
C^1(\overline{\Omega\backslash\omega})\cap 
C^1(\overline\omega)$.$\hfill\square$ 
\begin{thm} \label{lem:unique k +infty A} 
There exists at most one solution $u\in H^1(\Omega)\cap 
C^1(\overline{\Omega\backslash\omega})\cap C^1(\overline\omega)$ of 
equation (\ref{eq:k +infty A}). 
\end{thm} 
\emph{Proof}:\hspace{0.2cm}It is equivalent to showing that if 
$\varphi=0$, equation (\ref{eq:k +infty A}) only has the solution 
$u\equiv0$. Integrating by parts in the first line of equation 
(\ref{eq:k +infty A}), we have 
\begin{equation} 
\begin{aligned} 
0&=-\int_{\Omega\backslash\overline{\omega}}\partial_{x_j}\Big(a_2^{ij}(x)~\partial_{x_i}{u_k}\Big)\cdot u\\ 
&=\int_{\Omega\backslash\overline{\omega}}a_2^{ij}(x)\partial_{x_i}{u}\partial_{x_j}{u} 
- \int_{\partial\Omega}u\cdot 
a_2^{ij}(x)\partial_{x_i}{u}\nu_j\big|_{-} + 
\int_{\partial\omega}u\cdot a_2^{ij}(x)\partial_{x_i}{u}\nu_j\big|_{+}\\ 
&\geq\lambda\int_{\Omega\backslash\overline{\omega}}|\nabla u|^2 - 
\int_{\partial\Omega}\varphi\cdot 
a_2^{ij}(x)\partial_{x_i}{u}\nu_j\big|_{-}+ 
C_s\int_{\partial\omega_s}a_2^{ij}(x)\partial_{x_i}{u}\nu_j\big|_{+}\\ 
&=\lambda\int_{\Omega\backslash\overline{\omega}}|\nabla u|^2. 
\end{aligned} 
\nonumber 
\end{equation} 
Thus $\nabla u=0$ in $\Omega\backslash\overline{\omega}$. And since 
$u=\varphi=0$ on $\partial\Omega$, we have $u\equiv0$ in 
$\Omega\backslash\overline{\omega}$. Since $u|_{+}=u|_{-}$ on 
$\partial\omega$ and $u\equiv C$ on $\overline{\omega}$, we get 
$u=0$ on $\overline{\omega}$. Hence $u\equiv0$ in $\Omega$, i.e. 
$u\equiv0$ is the only solution of (\ref{eq:k +infty A}) when 
$\varphi=0$. $\hfill\square$\\ 
 
Define the energy functional 
\begin{equation} \label{vm:k +infty A} 
I_{\infty}[v]:= 
\frac{1}{2}\int_{\Omega\backslash\overline{\omega}}a_2^{ij}(x)\partial_{x_i}{v}\partial_{x_j}{v}, 
%\nonumber 
\end{equation} 
where $v$ belongs to the set 
$$ 
\mathcal{A}:=\big\{v\in H^1_\varphi(\Omega)\big| \nabla v\equiv 0 
~in~ \omega\big\}. 
$$ 
\begin{thm} \label{newlem:exist k finite A} 
There exists a minimizer $u\in \mathcal{A}$ satisfying 
$$I_{\infty}[u]=\min_{v\in \mathcal{A}}I_{\infty}[v].$$ 
Moreover, $u\in H^1(\Omega)\cap 
C^1(\overline{\Omega\backslash\omega})\cap C^1(\overline\omega)$ is 
a solution of equation (\ref{eq:k +infty A}). 
\end{thm} 
\emph{Proof}:\hspace{0.2cm} By the lower-semi continuity of 
$I_{\infty}$ and the weakly closed property of $\mathcal{A}$, it is 
easy to see that the minimizer $u\in\mathcal {A}$ exists and 
satisfies $\partial_{x_j}\Big(a_2^{ij}(x)\partial_{x_i}{u}\Big)=0$ 
in $\Omega\backslash\overline\omega$. The only thing which needs to 
be shown is the fourth line in equation (\ref{eq:k +infty A}), i.e. 
$$\int_{\partial\omega_s}a_2^{ij}(x)\partial_{x_i}{u}\nu_j\big|_{+}=0,~~~s=1,2,\cdots,m.$$ 
In fact, since $u$ is a minimizer, for any $\phi \in 
C_c^\infty(\Omega)$ satisfying $\phi\equiv 1$ on $\overline\omega_s$ 
and $\phi\equiv 0$ on $\overline\omega_t(t\neq s)$, let 
$$i(t):=I_{\infty}[u+t\phi]\hspace{0.5cm}(t\in\mathbb{R}),$$ 
we have 
$$i'(0):=\frac{di}{dt}\Big|_{t=0}=\int_{\Omega\backslash\overline{\omega}}a_2^{ij}(x)\partial_{x_i}{u}{\phi}_{x_j}=0.$$ 
Therefore 
\begin{equation} 
\begin{aligned} 
0&=-\int_{\Omega\backslash\overline{\omega}}\partial_{x_j}\Big(a_2^{ij}(x)~\partial_{x_i}{u_k}\Big)\phi 
=\int_{\Omega\backslash\overline{\omega}}a_2^{ij}(x)\partial_{x_i}{u}{\phi}_{x_j}+\int_{\partial\omega_s}\phi\cdot 
a_2^{ij}(x)\partial_{x_i}{u}\nu_j\big|_{+}\\ 
&=\int_{\partial\omega_s}a_2^{ij}(x)\partial_{x_i}{u}\nu_j\big|_{+}, 
\end{aligned} 
\nonumber 
\end{equation} 
for all $s=1,2,\cdots,m$.$\hfill\square$\\ 
 
Finally, we give the relationship between $u_k$ and $u$. 
\begin{thm} \label{thm:weak limit A} 
Let $u_k$ and $u$ in $H^1(\Omega)$ be the solutions of equations 
(\ref{eq:k finite A}) and (\ref{eq:k +infty A}), respectively. Then 
$$u_k\rightharpoonup u ~~\text{in}~H^1(\Omega), ~~~\text{as}~ k\rightarrow 
+\infty,$$ and 
$$\lim_{k\rightarrow +\infty}I_{k}[u_k]=I_{\infty}[u],$$ 
where $I_{k}$ and $I_{\infty}$ are defined as (\ref{vm:k finite A}) 
and (\ref{vm:k +infty A}). 
\end{thm} 
\emph{Proof}:\hspace{0.2cm} Step 1. By the uniqueness of the 
solution to equation (\ref{eq:k +infty A}), we only need to show 
that there exists a weak limit $u$ of a subsequence of $\{u_k\}$ in 
$H^1(\Omega)$ and $u$ is the solution of equation 
(\ref{eq:k +infty A}).\vspace{0.2cm}\\ 
(1) To show that after passing to a subsequence, $u_k$ weakly 
converges in $H^1(\Omega)$ to some $u$. 
 
\vspace{0.1cm}Let $\eta\in H^1_{\varphi}(\Omega)$ be fixed and 
satisfy $\eta\equiv 0$ on $\overline\omega$, then since $u_k$ is the 
minimizer of $I_{k}$ in $H^1_{\varphi}(\Omega)$, we have 
\begin{equation} 
\frac{\lambda}{2}\|\nabla u_k\|^2_{L^2(\Omega)} \leq I_{k}[u_k]\leq 
I_{k}[\eta]=\frac{1}{2}\int_{\Omega\backslash\overline\omega}a_2^{ij}(x)\eta_{x_i}\eta_{x_j} 
\leq\frac{\Lambda}{2}\|\eta\|^2_{H^1(\Omega)}, \nonumber 
\end{equation} 
i.e. 
$$\|\nabla u_k\|_{L^2(\Omega)}\leq\|\eta\|_{H^1(\Omega)}\doteq \overline{M},$$ 
where $\overline{M}$ is independent of $k$.\\ 
Since $u_k=\varphi$ on $\partial\Omega$ and 
$\sup_k\|u_k\|_{H^1(\Omega)}<\infty$, we have $u_k\rightharpoonup u 
\hspace{0.2cm}in\hspace{0.2cm}H^1_{\varphi}(\Omega).$\vspace{0.2cm}\\ 
(2) To show that $u$ is a solution of equation (\ref{eq:k +infty 
A}). 
 
\vspace{0.1cm}In fact, we only need to prove the following three 
conditions: 
\begin{align} 
\label{condition:u hamonic}\partial_{x_j}\Big(a_2^{ij}(x)~\partial_{x_i}{u}\Big)&=0~~~~~~~~in~\Omega\backslash\overline\omega,\\ 
\label{condition:u const}\nabla u&= 0~~~~~~~~in~\omega,\\ 
\label{condition:int 
0}\int_{\partial\omega_s}a_2^{ij}(x)\partial_{x_i}{u_k}\nu_j\big|_{+}&=0, 
~~~~~~s=1,2,\cdots,m. 
\end{align} 
(i) For every $k$, since $u_k\in H^1(\Omega)$ is the solution of 
equation (\ref{diveq:k finite A}), then 
 
$\forall\hspace{0.1cm}\phi\in C_c^\infty(\Omega),$ we have 
$$k\int_{\omega}a_1^{ij}(x)\partial_{x_i}{u_k}{\phi}_{x_j}+\int_{\Omega\backslash\overline{\omega}}a_2^{ij}(x)\partial_{x_i}{u_k}{\phi}_{x_j}=0.$$ 
Thus, $\forall\hspace{0.1cm}\phi\in 
C_c^\infty(\Omega\backslash\overline\omega)\subset 
C_c^\infty(\Omega),$ 
$$0=\int_{\Omega\backslash\overline{\omega}}a_2^{ij}(x)\partial_{x_i}{u_k}{\phi}_{x_j} 
\longrightarrow 
\int_{\Omega\backslash\overline{\omega}}a_2^{ij}(x)\partial_{x_i}{u}{\phi}_{x_j},$$ 
since $u_k\rightharpoonup u 
\hspace{0.2cm}in\hspace{0.2cm}H^1_{\varphi}(\Omega)\subset H^1(\Omega).$\\ 
Therefore, 
$$\int_{\Omega\backslash\overline{\omega}}a_2^{ij}(x)\partial_{x_i}{u}{\phi}_{x_j}=0,~~~\forall~\phi\in 
C_c^\infty(\Omega\backslash\overline\omega),$$ 
i.e. (\ref{condition:u hamonic}).\vspace{0.2cm}\\ 
(ii) Let $\eta\in H^1_{\varphi}(\Omega)$ be fixed and satisfy 
$\eta\equiv 0$ on $\overline\omega$, then since $u_k$ is the 
minimizer of $I_{k}$ in $H^1_{\varphi}(\Omega)$, we have 
\begin{equation} 
\frac{k\lambda}{2}\|\nabla u_k\|^2_{L^2(\omega)}\leq I_{k}[u_k]\leq 
I_{k}[\eta]=\frac{1}{2}\int_{\Omega\backslash\overline\omega}a_2^{ij}(x)\partial_{x_i}{\eta}\partial_{x_j}{\eta}\leq\frac{\Lambda}{2}\|\eta\|^2_{H^1(\Omega)}, 
\nonumber 
\end{equation} 
which implies 
$$\|\nabla u_k\|^2_{L^2(\omega)}\rightarrow 0,~~~~~~as ~k\rightarrow \infty.$$ 
By (1), since $u_k\rightharpoonup u ~in~H^1(\Omega),$ then 
$u_k\rightharpoonup u ~in~H^1(\omega).$ Therefore, by the lower-semi 
continuity, we get 
\begin{equation} 
\begin{split} 
0\leq\lambda\int_{\omega}|\nabla 
u|^2&\leq\int_{\omega}a_1^{ij}(x)\partial_{x_i}{u}\partial_{x_j}{u}\leq\int_{\omega}a_1^{ij}(x)\partial_{x_i}{u_k}\partial_{x_j}{u_k}\\ 
&\leq\Lambda\|\nabla u_k\|^2_{L^2(\omega)}\longrightarrow 0, 
\hspace{0.8cm}as\hspace{0.2cm}k\longrightarrow \infty . \nonumber 
\end{split} 
\end{equation} 
Hence, $~\int_{\omega}|\nabla u|^2=0 ~\Longrightarrow ~\nabla 
u\equiv 0$ in $\omega$, which is just 
(\ref{condition:u const}).\vspace{0.2cm}\\ 
(iii) By (i) and (ii), u satisfies (\ref{condition:u hamonic}) and 
is either constant or $\varphi$ on each component of 
$\partial(\Omega\backslash\overline\omega)$. Thus, $u\in 
C^{2}(\overline{\Omega\backslash\omega})$. For each 
$s=1,2,\cdots,m$, we construct a function $\varrho \in 
C^{2}(\overline{\Omega\backslash\omega})$, such that 
$\varrho=1$ on $\partial\omega_s$, $\varrho=0$ on $\partial\omega_t(t\neq s)$, and $\varrho=0$ on $\partial\Omega$.\\ 
By Green's Identity,  we have the following: 
\begin{equation} 
\begin{split} 
0&=-\int_{\Omega\backslash\overline{\omega}}\partial_{x_j}\Big(a_2^{ij}(x)~\partial_{x_i}{u_k}\Big)\varrho\\ 
&=\int_{\Omega\backslash\overline{\omega}}a_2^{ij}(x)\partial_{x_i}{u_k}\partial_{x_j}{\varrho}-\int_{\partial 
\Omega} \varrho\cdot 
a_2^{ij}(x)\partial_{x_i}{u_k}{\nu}_{j}\big|_{-}+\int_{\partial\omega}\varrho\cdot 
a_2^{ij}(x)\partial_{x_i}{u_k}{\nu}_{j}\big|_{+}\\ 
&=\int_{\Omega\backslash\overline{\omega}}a_2^{ij}(x)\partial_{x_i}{u_k}\partial_{x_j}{\varrho}+k\int_{\partial\omega_s} 
a_1^{ij}(x)\partial_{x_i}{u_k}{\nu}_{j}\big|_{-}\\ 
&=\int_{\Omega\backslash\overline{\omega}}a_2^{ij}(x)\partial_{x_i}{u_k}\partial_{x_j}{\varrho}. 
\nonumber 
\end{split} 
\end{equation} 
Similarly, 
$$0=-\int_{\Omega\backslash\overline{\omega}}\partial_{x_j}\Big(a_2^{ij}(x)~\partial_{x_i}{u}\Big)\varrho=\int_{\Omega\backslash\overline{\omega}} 
a_2^{ij}(x)\partial_{x_i}{u}\partial_{x_j}{\varrho}+\int_{\partial\omega_s} 
a_2^{ij}(x)\partial_{x_i}{u}{\nu}_{j}\big|_{+}.$$ Since 
$u_k\rightharpoonup u ~$in $H^1(\Omega)$, it 
follows$$0=\int_{\Omega\backslash\overline{\omega}}a_2^{ij}(x)\partial_{x_i}{u_k}\partial_{x_j}{\varrho}\longrightarrow\int_{\Omega\backslash\overline{\omega}}a_2^{ij}(x)\partial_{x_i}{u}\partial_{x_j}{\varrho}.$$ 
Thus, 
$$\int_{\partial\omega_s} 
a_2^{ij}(x)\partial_{x_i}{u}{\nu}_{j}\big|_{+}=0,$$ for any $s=1,2,\cdots,m$. Therefore, we finish 
the proof of the first part.\\ 
 
Step 2. Since $u_k$ is a minimizer of $I_{k}$ and $\nabla u=0$ in 
$\omega$, for any $k\in \mathbb{N}$, 
$$I_{k}[u_k]\leq I_{k}[u]=I_{\infty}[u].$$ 
Then $\limsup_{k\rightarrow +\infty}I_{k}[u_k]\leq I_{\infty}[u]$.\\ 
 
On the other hand, by Theorem \ref{thm:weak limit A}, since $u$ is 
the weak limit of $\{u_k\}$ in $H^1(\Omega)$, we obtain 
$$ 
I_{\infty}[u]=\int_{\Omega}a_2^{ij}(x)\partial_{x_i}{u}\partial_{x_j}{u}\leq 
\liminf_{k\rightarrow + 
\infty}\int_{\Omega}a_2^{ij}(x)\partial_{x_i}{u_k}\partial_{x_j}{u_k}\leq\liminf_{k\rightarrow 
+\infty}I_{k}[u_k]. 
$$ 
Therefore, 
$$\lim_{k\rightarrow +\infty}I_{k}[u_k]=I_{\infty}[u].$$ 
$\hfill\square$\\

\end{document}